\documentclass[preprint]{elsarticle}

\journal{arXiv}

\usepackage[dvipsnames]{xcolor}
\usepackage{amsmath,amsthm,amsfonts,amssymb,tikz}
\usepackage{hyperref}
\usepackage[capitalise,noabbrev]{cleveref}
\usepackage{multicol,multirow,diagbox,colortbl,arydshln}

\hypersetup{
    colorlinks,
    linkcolor={red!50!black},
    citecolor={red!50!black},
    urlcolor={blue!80!black}
}

\newcommand{\oper}[2]{\scriptstyle\binom{#1}{#2}_p}
\newcommand{\lc}[2]{\lambda^{#1}_{#2}}

\newcommand{\PG}[1]{{\foreach\x in #1{\x}}_p}
\newcommand{\PGs}[1]{\{#1\}_p}
\newcommand{\PN}[2]{#1^#2_p}
\newcommand{\PNN}[4]{{#1^#2 #3^#4}_p}
\newcommand{\lev}{\tl{\mbox{}}{\mbox{}}} 
\newcommand{\tl}[2]{\begin{tabular}[c]{@{}c@{}}#1\\#2 \end{tabular}}
\newcommand\E[1]{E_{#1}}

\DeclareMathOperator{\fitwtw}{\oper{5}{22}}
\DeclareMathOperator{\fifo}{\oper{5}{4}}
\DeclareMathOperator{\fotwon}{\oper{4}{21}}
\DeclareMathOperator{\foth}{\oper{4}{3}}
\DeclareMathOperator{\thtw}{\oper{3}{2}}
\DeclareMathOperator{\thonon}{\oper{3}{11}}
\DeclareMathOperator{\fithon}{\oper{5}{31}}
\DeclareMathOperator{\onze}{\oper{1}{0}}
\DeclareMathOperator{\twon}{\oper{2}{1}}
\DeclareMathOperator{\tmu}{\mu^{11}}
\DeclareMathOperator{\ms}{\mu^*}
\DeclareMathOperator{\les}{\le^*}
\DeclareMathOperator{\pf}{\mathcal{P}_{\le 5}^*}
\DeclareMathOperator{\pa}{\mathcal{P}^*}

\newtheorem{thm}{Theorem}\crefname{thm}{Theorem}{Theorems}
\newtheorem{lem}[thm]{Lemma}\crefname{lem}{Lemma}{Lemmas}
\newtheorem{ex}[thm]{Example}\crefname{ex}{Example}{Examples}
\newtheorem{cor}[thm]{Corollary}\crefname{cor}{Corollary}{Corollaries}
\newtheorem{conj}[thm]{Conjecture}\crefname{conj}{Conjecture}{Conjectures}
\newtheorem{prop}[thm]{Proposition}\crefname{prop}{Proposition}{Propositions}
\newtheorem{defn}[thm]{Definition}\crefname{defn}{Definition}{Definitions}
\newtheorem{rem}[thm]{Remark}\crefname{rem}{Remark}{Remarks}
\newtheorem{que}[thm]{Question}\crefname{que}{Question}{Questions}
\newtheorem{clm}{Claim}\crefname{clm}{Claim}{Claims}
\crefname{figure}{Figure}{Figure}

\numberwithin{equation}{section}
\numberwithin{figure}{section}
\numberwithin{thm}{section}

\begin{document}

\begin{frontmatter}
\title{The poset of graphs ordered by induced containment}

\author{Jason P. Smith\fnref{fn1}}
\address{Department of Computer and Information Sciences, University of Strathclyde, \\
             Glasgow, UK.}
\ead{jason.p.smith@strath.ac.uk}
\fntext[fn1]{This research was supported by the EPSRC Grant EP/M027147/1.}

\begin{abstract}
We study the poset $\mathcal{G}$ of all unlabelled graphs, up to isomorphism, with $H\le G$ if $H$ occurs
 as an induced subgraph in $G$. We present some general results on the M\"obius function
 of intervals of $\mathcal{G}$ and some results for specific classes of graphs. This includes a case where the 
 M\"obius function is given by the Catalan numbers, which we prove using discrete Morse theory, and another case
 where it equals the Fibonacci numbers, therefore showing that the M\"obius 
 function is unbounded. A classification of the disconnected intervals of $\mathcal{G}$ is presented,
 which gives a large class of non-shellable intervals. We also present several conjectures on the structure of $\mathcal{G}$.
\end{abstract}
\end{frontmatter}

\section{Introduction}
Given any set of combinatorial objects and a suitable notion of containment of one such object in another we can define a poset. Many such posets have been
 studied in the literature, such as the posets of words with subword order~\cite{Bjo90}, the permutation pattern poset \cite{BJJS11,Smith15},
 the poset of graphs with minor order \cite{RobSey04}, and many more. In this paper we introduce and study the poset $\mathcal{G}$ of all graphs ordered by 
 induced containment, that is, $\mathcal{G}$ contains all unlabelled finite graphs and $H\le G$ if $H$ is an induced 
 subgraph of~$G$. We say that $G$ contains $H$ if $H\le G$ and we consider two subgraphs to be the same if they are 
 isomorphic.

We study the M\"obius function and topology of $\mathcal{G}$, and present some results and conjectures.  
 The poset $\mathcal{G}$ has a countably infinite number of elements and is locally finite, so we focus our attention 
 on the \emph{intervals}~$[a,b]={\{z\in\mathcal{G}\,|\,a\le z\le b\}}$, see \cref{fig:2C4}. We also consider the poset~$\mathcal{G}^c$
 of all connected graphs, which is an induced subposet of $\mathcal{G}$, and we denote intervals of $\mathcal{G}^c$ by~$[a,b]^c$.
 Unless otherwise specified we allow graphs to have loops and multiple edges.

Whilst the topology of the poset $\mathcal{G}$ does not seem to have been studied, the poset has been considered from a model theoretic perspective in \cite{Wir16}. There are also many other posets of graphs that have been investigated.
 In \cite{Tha06} the poset of all induced connected subgraphs of a graph~$G$ is considered in relation to the graph reconstruction conjecture.
 This is equivalent to the interval~$[K_1,G]^c$ in~$\mathcal{G}^c$, where~$K_1$ is the graph with a single vertex, 
 but the topology and M\"obius function of these posets is not considered. In \cite{Kez96} the poset $C(G)$ is defined 
 on a graph $G$ with the same partial order considered here, but graphs are not considered the same if they are isomorphic,
 such a poset has a simpler structure than $\mathcal{G}$ because subgraphs can occur exactly once in the parent graph.
 In \cite{Smi17} the posets of connected labelled graphs
 on $n$ vertices ordered by non-induced containment is shown to be Sperner. 
 Another poset that has been considered previously is that of graphs ordered by the minor relation, that is, $G\le H$ if~$G$
 is a graphical minor of $H$. This is well known to be a partial order and was famously shown to be a well-quasi
 ordering by the Robertson-Seymour Theorem~\cite{RobSey04}. 

A finite poset $P$ is \emph{shellable} if there is a ``nice" ordering of the maximal \emph{chains}, that is, the maximal totally 
 ordered subsets of~$P$, see \cite{Wac07} for a formal definition of shellability. A poset is not shellable if it contains 
 any disconnected subintervals of rank at least $3$. We give a classification of the disconnected intervals of $\mathcal{G}$, 
 which is similar to that given for permutation patterns in~\cite{McSt13}. Moreover, we present a large class of non-shellable 
 intervals of $\mathcal{G}$, but we conjecture that the proportion of intervals that are shellable tends to $1$ as the rank of 
 the intervals increases.
 
 Note that $\mathcal{G}$ is a ranked poset, that is, all maximal chain have the same length, and is not a lattice. It is straightforward to see that the rank function is
 simply the order $|G|$ of a graph $G$, that is, the number of vertices. Moreover, it can easily be seen in \cref{fig:2C4} that $\mathcal{G}$
 is not a lattice, as there are multiple pairs which do not have a unique join.

In \cref{sec:conn} we give a classification of the disconnected intervals $[H,G]$ of $\mathcal{G}$ based on the set of occurrences
 of $H$ in $G$. In \cref{sec:MF} we give some general results on the M\"obius function of intervals of~$\mathcal{G}$.
 In \cref{sec:known} 
 we consider intervals of $\mathcal{G}$ and $\mathcal{G}^c$ between some well known graphs, such as the complete graphs $K_n$, 
 the cycle graphs $C_n$ and the empty graphs $\overline{K}_n$, that is, the graphs with no edges. In 
 \cref{sec:paths} we consider intervals between graphs of disjoint paths. We prove a result where the M\"obius function 
 is given by the Catalan numbers, using discrete Morse theory, and another case where it is given by the Fibonacci numbers.
 As a corollary we get that the M\"obius function is unbounded on $\mathcal{G}$. 
 In \cref{sec:conj} we finish with some conjectures about $\mathcal{G}$.

\section{Disconnected Intervals}\label{sec:conn}
In this section we consider the disconnected intervals of $\mathcal{G}$ and $\mathcal{G}^c$.
 An interval~$[x,y]$ is \emph{disconnected} if the interior $(x,y):=[x,y]\setminus\{x,y\}$
 can be split into two non-empty sets $A$ and $B$, 
 which we call $\emph{components}$, with $a\not\le b$ and~$b\not\le a$ for all~$a\in A$ and $b\in B$. A disconnected
 subinterval is \emph{non-trivial} if it has a rank of at least $3$, and it is shown in \cite[Proposition~4.2]{Bjo80} that 
 a non-trivial disconnected poset is not shellable.

In \cite{McSt13} a classification of the disconnected intervals of the permutation poset is given based on
 splitting the set of occurrences into two disjoint sets, and this result is generalised in \cite{Smith16} to general 
 pattern posets. In this section we introduce and apply an analogous result to the poset of graphs. First we need to define an 
 occurrence in this setting.

Given a graph $G$ we denote the set of vertices and set of edges of $G$ by $V(G)$ and $E(G)$, respectively.
 We arbitrarily assign the labels $[|G|]:=\{1,\ldots,|G|\}$ to the elements of $V(G)$. This labelling is just to 
 record which vertices give an induced subgraph $H$, but we do not consider the graph as a labelled graph. Given a 
 set $\eta\subseteq [|G|]$ we denote the subgraph of $G$ induced by $\eta$ as $G[\eta]$ and say that $\eta$ is an 
 \emph{occurrence} of $H$ in $G$ if $G[\eta]\cong H$. Let $\E{H,G}$ be the set of all occurrences of $H$ in $G$. 
 Given any $\eta\in \E{H,G }$ define the set $Z(\eta):=[|G|]\setminus\eta$ and for any $A\subseteq \E{H,G}$ let 
 $Z(A)=\bigcup_{\eta\in A}Z(\eta)$. For ease of notation we will often denote an occurrence $\{a,b,c\}$ 
 as simply $abc$. Next we introduce the notion \emph{zero-splitness}, first introduced in \cite[Proposition 5.3]{McSt13}, although 
 not with the same name.

\begin{figure}\centering
\def\s{0.4}
\begin{tikzpicture}[scale=2.1]\def\y{1.5}
\node (41) at (0,4*\y){\begin{tikzpicture}[scale=0.5]
\node[circle, fill=black,scale=\s] (1) at (0,0){};
\node[circle, fill=black,scale=\s] (3) at (0,1){};
\node[circle, fill=black,scale=\s] (2) at (1,0){};
\node[circle, fill=black,scale=\s] (4) at (1,1){};
\node[circle, fill=black,scale=\s] (5) at (2,0.5){};
\node[circle, fill=black,scale=\s] (6) at (2,1.5){};
\node[circle, fill=black,scale=\s] (7) at (3,0.5){};
\node[circle, fill=black,scale=\s] (8) at (3,1.5){};
\draw[thick] (8) -- (2) -- (1) -- (3) -- (4) -- (2) -- (5) -- (6) -- (8) -- (7) -- (5);
\draw[thick] (7) -- (2) -- (6);
\end{tikzpicture}};
\node (31) at (-1.5,3*\y){\begin{tikzpicture}[scale=0.5]
\node[circle, fill=black,scale=\s] (1) at (0,0){};
\node[circle, fill=black,scale=\s] (3) at (0,1){};
\node[circle, fill=black,scale=\s] (2) at (1,0){};
\node[circle, fill=black,scale=\s] (4) at (1,1){};
\node[circle, fill=black,scale=\s] (5) at (2,0.5){};
\node[circle, fill=black,scale=\s] (6) at (2,1.5){};
\node[circle, fill=black,scale=\s] (7) at (3,0.5){};
\draw[thick] (2) -- (1) -- (3) -- (4) -- (2) -- (5) -- (6);
\draw[thick] (5) -- (7) -- (2) -- (6);
\end{tikzpicture}};
\node (32) at (-0.5,3*\y){\begin{tikzpicture}[scale=0.5]
\node[circle, fill=black,scale=\s] (3) at (0,1){};
\node[circle, fill=black,scale=\s] (2) at (1,0){};
\node[circle, fill=black,scale=\s] (4) at (1,1){};
\node[circle, fill=black,scale=\s] (5) at (2,0.5){};
\node[circle, fill=black,scale=\s] (6) at (2,1.5){};
\node[circle, fill=black,scale=\s] (7) at (3,0.5){};
\node[circle, fill=black,scale=\s] (8) at (3,1.5){};
\draw[thick] (8) -- (2);\draw[thick] (3) -- (4) -- (2) -- (5) -- (6) -- (8) -- (7) -- (5);
\draw[thick] (7) -- (2) -- (6);
\end{tikzpicture}};
\node (33) at (0.5,3*\y){\begin{tikzpicture}[scale=0.5]
\node[circle, fill=black,scale=\s] (1) at (0,0){};
\node[circle, fill=black,scale=\s] (2) at (1,0){};
\node[circle, fill=black,scale=\s] (4) at (1,1){};
\node[circle, fill=black,scale=\s] (5) at (2,0.5){};
\node[circle, fill=black,scale=\s] (6) at (2,1.5){};
\node[circle, fill=black,scale=\s] (7) at (3,0.5){};
\node[circle, fill=black,scale=\s] (8) at (3,1.5){};
\draw[thick] (8) -- (2) -- (1);\draw[thick]  (4) -- (2) -- (5) -- (6) -- (8) -- (7) -- (5);
\draw[thick] (7) -- (2) -- (6);
\end{tikzpicture}};
\node (34) at (1.5,3*\y){\begin{tikzpicture}[scale=0.5]
\node[circle, fill=black,scale=\s] (1) at (0,0){};
\node[circle, fill=black,scale=\s] (3) at (0,1){};
\node[circle, fill=black,scale=\s] (4) at (1,1){};
\node[circle, fill=black,scale=\s] (5) at (2,0.5){};
\node[circle, fill=black,scale=\s] (6) at (2,1.5){};
\node[circle, fill=black,scale=\s] (7) at (3,0.5){};
\node[circle, fill=black,scale=\s] (8) at (3,1.5){};
\draw[thick] (1) -- (3) -- (4);\draw[thick] (5) -- (6) -- (8) -- (7) -- (5);
\end{tikzpicture}};
\node (21) at (-2.5,2*\y){\begin{tikzpicture}[scale=0.5]
\node[circle, fill=black,scale=\s] (1) at (0,0){};
\node[circle, fill=black,scale=\s] (3) at (0,1){};
\node[circle, fill=black,scale=\s] (2) at (1,0){};
\node[circle, fill=black,scale=\s] (4) at (1,1){};
\node[circle, fill=black,scale=\s] (5) at (2,0.5){};
\node[circle, fill=black,scale=\s] (6) at (2,1.5){};
\draw[thick] (6) -- (2) -- (1) -- (3) -- (4) -- (2) -- (5) -- (6);
\end{tikzpicture}};
\node (22) at (-1.5,2*\y){\begin{tikzpicture}[scale=0.5]
\node[circle, fill=black,scale=\s] (1) at (0,0){};
\node[circle, fill=black,scale=\s] (3) at (0,1){};
\node[circle, fill=black,scale=\s] (2) at (1,0){};
\node[circle, fill=black,scale=\s] (4) at (1,1){};
\node[circle, fill=black,scale=\s] (6) at (2,1.5){};
\node[circle, fill=black,scale=\s] (7) at (3,0.5){};
\draw[thick] (2) -- (1) -- (3) -- (4) -- (2);
\draw[thick] (7) -- (2) -- (6);
\end{tikzpicture}};
\node (23) at (-0.5,2*\y){\begin{tikzpicture}[scale=0.5]
\node[circle, fill=black,scale=\s] (3) at (1,1){};
\node[circle, fill=black,scale=\s] (2) at (1,0){};
\node[circle, fill=black,scale=\s] (5) at (2,0.5){};
\node[circle, fill=black,scale=\s] (6) at (2,1.5){};
\node[circle, fill=black,scale=\s] (7) at (3,0.5){};
\node[circle, fill=black,scale=\s] (8) at (3,1.5){};
\draw[thick] (8) -- (2) -- (5) -- (6) -- (8) -- (7) -- (5);
\draw[thick] (7) -- (2) -- (6);
\end{tikzpicture}};
\node (24) at (0.5,2*\y){\begin{tikzpicture}[scale=0.5]
\node[circle, fill=black,scale=\s] (2) at (1,0){};
\node[circle, fill=black,scale=\s] (4) at (1,1){};
\node[circle, fill=black,scale=\s] (5) at (2,0.5){};
\node[circle, fill=black,scale=\s] (6) at (2,1.5){};
\node[circle, fill=black,scale=\s] (7) at (3,0.5){};
\node[circle, fill=black,scale=\s] (8) at (3,1.5){};
\draw[thick] (8) -- (2);\draw[thick] (4) -- (2) -- (5) -- (6) -- (8) -- (7) -- (5);
\draw[thick] (7) -- (2) -- (6);
\end{tikzpicture}};
\node (25) at (1.5,2*\y){\begin{tikzpicture}[scale=0.5]
\node[circle, fill=black,scale=\s] (3) at (0,1){};
\node[circle, fill=black,scale=\s] (4) at (1,1){};
\node[circle, fill=black,scale=\s] (5) at (2,0.5){};
\node[circle, fill=black,scale=\s] (6) at (2,1.5){};
\node[circle, fill=black,scale=\s] (7) at (3,0.5){};
\node[circle, fill=black,scale=\s] (8) at (3,1.5){};
\draw[thick] (3) -- (4);
\draw[thick] (5) -- (6) -- (8) -- (7) -- (5);
\end{tikzpicture}};
\node (26) at (2.5,2*\y){\begin{tikzpicture}[scale=0.5]
\node[circle, fill=black,scale=\s] (1) at (0,0){};
\node[circle, fill=black,scale=\s] (4) at (1,1){};
\node[circle, fill=black,scale=\s] (5) at (2,0.5){};
\node[circle, fill=black,scale=\s] (6) at (2,1.5){};
\node[circle, fill=black,scale=\s] (7) at (3,0.5){};
\node[circle, fill=black,scale=\s] (8) at (3,1.5){};
\draw[thick] (5) -- (6) -- (8) -- (7) -- (5);
\end{tikzpicture}};
\node (11) at (-1,1*\y){\begin{tikzpicture}[scale=0.5]
\node[circle, fill=black,scale=\s] (1) at (0,0){};
\node[circle, fill=black,scale=\s] (3) at (0,1){};
\node[circle, fill=black,scale=\s] (2) at (1,0){};
\node[circle, fill=black,scale=\s] (4) at (1,1){};
\node[circle, fill=black,scale=\s] (5) at (2,0.5){};
\draw[thick] (2) -- (1) -- (3) -- (4) -- (2) -- (5);
\end{tikzpicture}};
\node (12) at (0,1*\y){\begin{tikzpicture}[scale=0.5]
\node[circle, fill=black,scale=\s] (2) at (1,0){};
\node[circle, fill=black,scale=\s] (5) at (2,0.5){};
\node[circle, fill=black,scale=\s] (6) at (2,1.5){};
\node[circle, fill=black,scale=\s] (7) at (3,0.5){};
\node[circle, fill=black,scale=\s] (8) at (3,1.5){};
\draw[thick] (8) -- (2) -- (5) -- (6) -- (8) -- (7) -- (5);
\draw[thick] (7) -- (2) -- (6);
\end{tikzpicture}};
\node (13) at (1,1*\y){\begin{tikzpicture}[scale=0.5]
\node[circle, fill=black,scale=\s] (4) at (1,1){};
\node[circle, fill=black,scale=\s] (5) at (2,0.5){};
\node[circle, fill=black,scale=\s] (6) at (2,1.5){};
\node[circle, fill=black,scale=\s] (7) at (3,0.5){};
\node[circle, fill=black,scale=\s] (8) at (3,1.5){};
\draw[thick] (5) -- (6) -- (8) -- (7) -- (5);
\end{tikzpicture}};
\node (01) at (0,0*\y){\begin{tikzpicture}[scale=0.5]
\node[circle, fill=black,scale=0.5] (1) at (0,0){};
\node[circle, fill=black,scale=0.5] (3) at (0,1){};
\node[circle, fill=black,scale=0.5] (2) at (1,0){};
\node[circle, fill=black,scale=0.5] (4) at (1,1){};
\draw[thick] (2) -- (1) -- (3) -- (4) -- (2);
\end{tikzpicture}};
\draw[color=Gray] (41) -- (31) -- (21) -- (11) -- (22) -- (31);
\draw[color=Gray] (11) -- (01) -- (12) -- (23) -- (32) -- (41) -- (33) -- (24) -- (12);
\draw[color=Gray] (01) -- (13) -- (25) -- (34) -- (41);
\draw[color=Gray] (23) -- (13) -- (24) -- (32) -- (25);
\draw[color=Gray] (13) -- (26) -- (34);
\draw[color=Gray] (33) -- (26);
\end{tikzpicture}
\caption{The disconnected interval $[C_4,H]$ in $\mathcal{G}$,  for any $v\in V(C_4)$, where $H=D_v(C_4)$ as defined in \cref{def:D}.}\label{fig:2C4}
\end{figure}
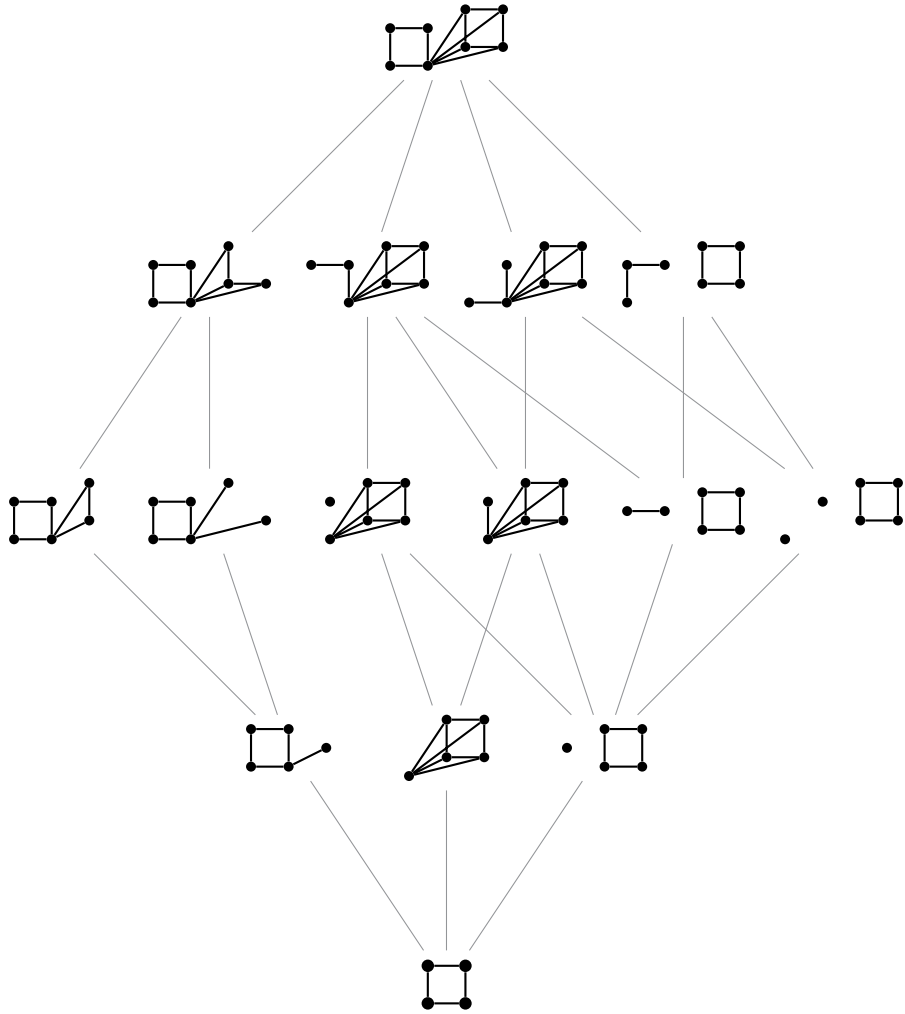

\begin{defn}
The interval $[H,G]$ is \emph{zero-split} if $\E{H,G}$ can be split into two non-empty disjoint sets $A$ and $B$ such
 that $Z(A)\cap Z(B)=\emptyset$. We call $A$ and $B$ a \emph{zero-split partition} of $[H,G]$.

The interval $[H,G]$ is \emph{strongly zero-split} if there is a zero-split partition~$A$ and $B$ and there
 does not exist $\eta\in A$, $\phi\in B$ and $i,j\in[|G|]$ such that ${G[\eta\cup i]\cong G[\phi\cup j]}$.
\end{defn}

\begin{ex}
Consider the interval $[C_4,G]$ where
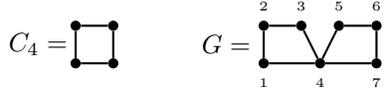
\begin{center}
\begin{tikzpicture}[scale=0.5]\def\s{0.4}\def\x{-5}
\node[circle, fill=black,scale=\s] (a) at (0+\x,0){};
\node[circle, fill=black,scale=\s] (b) at (0+\x,1){};
\node[circle, fill=black,scale=\s] (d) at (1+\x,0){};
\node[circle, fill=black,scale=\s] (c) at (1+\x,1){};
\draw[thick] (a) -- (b) -- (c) -- (d) -- (a);
\node (lab) at (-1+\x,0.5){$C_4=$};
\node[label={[font=\tiny]below:$1$},circle, fill=black,scale=\s] (1) at (0,0){};
\node[label={[font=\tiny]above:$2$},circle, fill=black,scale=\s] (3) at (0,1){};
\node[label={[font=\tiny]below:$4$},circle, fill=black,scale=\s] (2) at (1.5,0){};
\node[label={[font=\tiny]above:$3$},circle, fill=black,scale=\s] (4) at (1,1){};
\node[label={[font=\tiny]above:$5$},circle, fill=black,scale=\s] (6) at (2,1){};
\node[label={[font=\tiny]below:$7$},circle, fill=black,scale=\s] (7) at (3,0){};
\node[label={[font=\tiny]above:$6$},circle, fill=black,scale=\s] (8) at (3,1){};
\draw[thick] (2) -- (1) -- (3) -- (4) -- (2) -- (6) -- (8) -- (7) -- (2);
\node (lab) at (-1,0.5){$G=$};
\end{tikzpicture}.
\end{center}
So $\E{C_4,G}=\{1234,4567\}$ and there is a zero-split partition $A=\{1234\}$ and ${B=\{4567\}}$.
 However, the interval is not strongly zero-split because
 $$G[1234\cup7]\cong G[4567\cup1]\cong\raisebox{-5pt}{\begin{tikzpicture}[scale=0.5]\def\s{0.4}
  \node[circle, fill=black,scale=\s] (a) at (0,0){};
  \node[circle, fill=black,scale=\s] (b) at (0,1){};
  \node[circle, fill=black,scale=\s] (d) at (1,0){};
  \node[circle, fill=black,scale=\s] (c) at (1,1){};
  \node[circle, fill=black,scale=\s] (e) at (2,0){};
  \draw[thick] (a) -- (b) -- (c) -- (d) -- (a);\draw[thick] (d) -- (e);
 \end{tikzpicture}}.$$
\end{ex}

We can now present the main theorem of this section.

\begin{thm}\label{lem:discon}
An interval $[H,G]$, with $|G|-|H|>2$, is
 disconnected if and only if it is strongly zero-split.
\begin{proof}
Suppose $[H,G]$ is strongly zero-split with partition $E_1$ and $E_2$.
 Define the sets $P_1$ and $P_2$ such that $X\in P_i$ if there is an $\alpha\in \E{X,G}$
  with $Z(\alpha)\subset Z(\eta)$ for some $\eta\in E_i$. Every graph in $(H,G)$ can be
  obtained by deleting some subset of $Z(\eta)$ for some $\eta\in \E{H,G}$, so~$(H,G)=P_1\cup P_2$.

 We claim that $P_1$ and $P_2$ are disconnected components of $(H,G)$. First note
  that if $X\in P_i$ and $Y\le X$, then $Y\in P_i$. To see this consider $\alpha\in \E{X,G}$ 
  with $Z(\alpha)\subset Z(\eta)$ for some~$\eta\in E_i$. There is a $\beta\in \E{Y,G}$ with 
  $Z(\beta)\supseteq Z(\alpha)$. Moreover, $Z(\beta)\subset Z(\zeta)$ for some $\zeta\in \E{H,G}$ 
  which implies $Z(\zeta)\cap Z(\eta)\supseteq Z(\alpha)$. As $Z(\eta)$ and $Z(\zeta)$ have non-empty 
  intersection, $\eta$ and $\zeta$ cannot be in separate sets of the partition, which implies 
  $\zeta\in E_i$ thus $Y\in P_i$.

For a contradiction suppose $P_1$ and $P_2$ are not disconnected so there exists a comparable 
 pair~$p_1\in P_1$ and $p_2\in P_2$, and without loss of generality suppose~$p_1\le p_2$. This implies that 
 $p_1\in P_1\cap P_2$, and consider any $A\le p_1$ with~$|A|=|H|+1$, then $A\in P_1\cap P_2$. However, 
 this violates the strongly zero-split condition as it implies there is a pair $e_1\in E_1$ and $e_2\in E_2$ 
 both of which are contained in occurrences of $A$ in $G$. Therefore, $P_1$ and $P_2$ are disconnected components.

Suppose that $[H,G]$ is disconnected with components $P_1$ and $P_2$, and let 
$$E_i=\left\{\eta\in \E{H,G}\,\,\middle|\,\,Z(\eta)\subseteq \bigcup_{B\in P_i}Z\left(\E{B,G}\right)\right\}.$$
So $E_i$ is the set of occurrences $\eta$ from which an element of $P_i$ can be obtained by removing a subset of 
$Z(\eta)$ from $G$. The sets $E_1$ and $E_2$ form a zero-split partition because if $i\in Z(E_1)\cap Z(E_2)$ then 
$G-i$ is in both $P_1$ and $P_2$, which gives a contradiction. To see this is a strongly zero-split partition 
suppose there is a graph $X$ that can be obtained by adding a vertex to occurrences from different sets $E_1$ 
and $E_2$, this would imply that $X\in P_1\cap P_2$, again giving a contradiction. Therefore, $[H,G]$ is strongly zero-split.
\end{proof}
\end{thm}

\begin{cor}
If $|G|>2|H|$ then $[H,G]$ is connected.
\begin{proof}
For every pair $\eta,\zeta\in \E{H,G}$ we have $|Z(\eta)\cap Z(\zeta)|>0$ so $[H,G]$ is not strongly 
zero-split and the result follows by \cref{lem:discon}.
\end{proof}
\end{cor}

We get the following corollary by applying an analogous proof to that used in \cref{lem:discon}:
\begin{cor}
Consider an interval $[H,G]^c$ in $\mathcal{G}^c$ with $|G|-|H|>2$, then $[H,G]^c$ is disconnected if 
and only if it is strongly zero-split.
\end{cor}

We can use \cref{lem:discon} to construct an infinite class of disconnected intervals, using the following construction.
\begin{defn}\label{def:D}
Given a graph $H$ 
 and vertex $v\in V(H)$ define the graph $D_v(H)$ with 
 vertex and edge sets: 
 $$V(D_v(H))=\{1,2\}\times V(H),$$
 \begin{align*}E(D_v(H))=\{(i,x),(i,y)\,|\,(x,y)\in &E(H), i\in \{1,2\}\}\\&\cup\{((1,v),(2,x))\,|\,x\in V(H)\},\end{align*}
 that is, $D_v(H)$ is constructed by taking two copies of $H$ and connecting all vertices of one copy to the 
 vertex $v$ in the other copy.
\end{defn} 
See \cref{fig:2C4} for an example of $D_v(H)$. Define a \emph{pendant} in a 
 graph as a vertex with exactly one neighbour.

\begin{lem}\label{lem:disEx}
Consider a connected graph $H$, with $|H|\ge3$, that does not contain~$C_3$ and has no pendants. 
 The interval $[H,D_v(H)]$ is disconnected, for all~${v\in V(H)}$.
\begin{proof}
There are clearly two occurrences of $H$ in $D_v(H)$, namely the two copies of $H$ that $D_v(H)$
 is constructed from. Let $A$ be the occurrence of $H$ containing~$v$, and $B$ be the other occurrence.
 Suppose there is a third occurrence~$C$, then $C$ must contain vertices of both $A$ and $B$,
 and must contain $v$, since otherwise it is disconnected.
 If there are no edges between elements of $C\cap B$, then each of these vertices has a single
 neighbour in $C$ which is $v$, thus they are pendants in~$H$ which is not allowed.
 If there are edges between two elements $x,y\in C\cap B$, then 
 $x,y,v$ is an occurrence of $C_3$ in $C$, and thus in $H$ which is not allowed. Therefore, there is no third occurrence of $H$.

 By the definition of $D_v(H)$ the occurrences $A$ and $B$
 are zero-split. Moreover, in $A\cup i$ the vertex $i$ is a pendant for any $i\in D_v(H)\setminus A$,
 in $B\cup j$ the vertex~$j$ has $|H|$ neighbours if $j=v$ and no neighbours for any
 other $j\in D_v(H)\setminus B$. Therefore $A\cup i\not\cong B\cup j$ for any $i$ and $j$, so
 $[H,D_v(H)]$ is strongly zero-split, and thus disconnected by \cref{lem:discon}.
\end{proof}
\end{lem}

See \cref{fig:2C4} for an example of \cref{lem:disEx}. A poset is non-shellable if it contains a
 non-trivial disconnected subinterval, that is, a subinterval of rank at least $3$. Therefore, we can 
 use \cref{lem:disEx} to get an infinite class of non-shellable intervals.

\begin{cor}\label{cor:disEx}
Consider a connected graph $H$, with $|H|\ge3$, that does not contain $C_3$ and has no pendants. As $n$
 tends to infinity the probability that~$[H,G]$, where $|G|=n$, contains a non-trivial disconnected 
 subinterval, and thus is not shellable, tends to $1$.
\begin{proof}
For any fixed graph $H$ the probability that $H$ occurs as an induced subgraph in a graph $G$ tends to $1$ as $|G|$
 tends to infinity, see \cite[Section 11, Exercise~12]{Die12}.
 Therefore, the probability that $G$ contains $D_v(H)$ tends to $1$ as~$|G|$ tends to infinity, 
 which implies $[H,G]$ contains the subinterval $[H,D_v(H)]$, which is disconnected by \cref{lem:disEx}.
\end{proof}
\end{cor}

\section{M\"obius Function}\label{sec:MF}
In this section we give some results on the M\"obius function of intervals of~$\mathcal{G}$.
 Define the half-open interval $[a,b):=[a,b]\setminus\{b\}$.
 The \emph{M\"obius function} of a poset~$P$ is defined recursively by $\mu_P(a,a)=1$ for all~$a$, 
 $\mu_P(a,b)=0$ if $a\not\le b$ and if $a<b$ then 
 $$\mu_P(a,b)=-\sum_{c\in[a,b)}\mu_P(a,c).$$
 
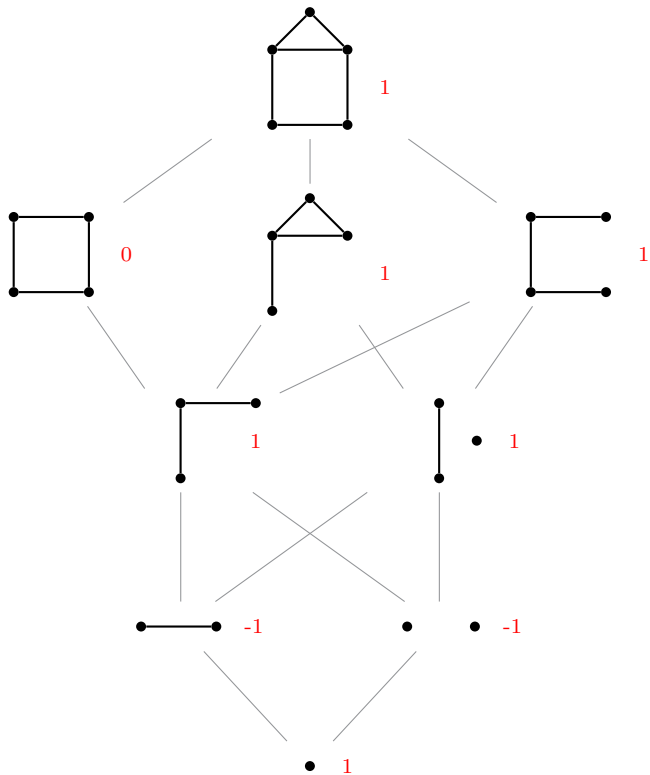
\begin{figure}[ht]\centering
\begin{tikzpicture}[scale=2.75]\def\s{0.4}\def\x{1.25}\def\y{0.9}
\node (H) at (0*\x,4*\y){\begin{tikzpicture}
\node[circle, fill=black,scale=\s] (1) at (0,3.5){};
\node[circle, fill=black,scale=\s] (2) at (-.5,3){};
\node[circle, fill=black,scale=\s] (3) at (.5,3){};
\node[circle, fill=black,scale=\s] (4) at (-.5,2){};
\node[circle, fill=black,scale=\s] (5) at (.5,2){};
\draw[thick] (2) -- (3) -- (1) -- (2) -- (4) -- (5) -- (3);
\node (lab) at (1,2.5){\footnotesize\textcolor{red}{1}};
\node (labi) at (-1,2.5){\footnotesize\textcolor{white}{0}};
\end{tikzpicture}};
\node (C4) at (-1*\x,3*\y){\begin{tikzpicture}
\node[circle, fill=black,scale=\s] (2) at (-.5,3){};
\node[circle, fill=black,scale=\s] (3) at (.5,3){};
\node[circle, fill=black,scale=\s] (4) at (-.5,2){};
\node[circle, fill=black,scale=\s] (5) at (.5,2){};
\draw[thick] (3) -- (2) -- (4) -- (5) -- (3);
\node (lab) at (1,2.5){\footnotesize\textcolor{red}{0}};
\node (labi) at (-1,2.5){\footnotesize\textcolor{white}{0}};
\end{tikzpicture}};
\node (F) at (0*\x,3*\y){\begin{tikzpicture}
\node[circle, fill=black,scale=\s] (1) at (0,3.5){};
\node[circle, fill=black,scale=\s] (2) at (-.5,3){};
\node[circle, fill=black,scale=\s] (3) at (.5,3){};
\node[circle, fill=black,scale=\s] (4) at (-.5,2){};
\draw[thick] (2) -- (3) -- (1) -- (2) -- (4);
\node (lab) at (1,2.5){\footnotesize\textcolor{red}{1}};
\node (labi) at (-1,2.5){\footnotesize\textcolor{white}{0}};
\end{tikzpicture}};
\node (P4) at (1*\x,3*\y){\begin{tikzpicture}
\node[circle, fill=black,scale=\s] (2) at (-.5,3){};
\node[circle, fill=black,scale=\s] (3) at (.5,3){};
\node[circle, fill=black,scale=\s] (4) at (-.5,2){};
\node[circle, fill=black,scale=\s] (5) at (.5,2){};
\draw[thick] (3) -- (2) -- (4) -- (5);
\node (lab) at (1,2.5){\footnotesize\textcolor{red}{1}};
\node (labi) at (-1,2.5){\footnotesize\textcolor{white}{0}};
\end{tikzpicture}};
\node (P3) at (-0.5*\x,2*\y){\begin{tikzpicture}
\node[circle, fill=black,scale=\s] (2) at (1,1){};
\node[circle, fill=black,scale=\s] (1) at (0,1){};
\node[circle, fill=black,scale=\s] (0) at (0,0){};
\draw[thick] (0) -- (1) -- (2);
\node (lab) at (1,.5){\footnotesize\textcolor{red}{1}};
\node (labi) at (-1,.5){\footnotesize\textcolor{white}{0}};
\end{tikzpicture}};
\node (21) at (0.5*\x,2*\y){\begin{tikzpicture}
\node[circle, fill=black,scale=\s] (2) at (0.5,0.5){};
\node[circle, fill=black,scale=\s] (1) at (0,1){};
\node[circle, fill=black,scale=\s] (0) at (0,0){};
\draw[thick] (0) -- (1);
\node (lab) at (1,.5){\footnotesize\textcolor{red}{1}};
\node (labi) at (-1,.5){\footnotesize\textcolor{white}{0}};
\end{tikzpicture}};
\node (P2) at (-0.5*\x,1*\y){\begin{tikzpicture}
\node[circle, fill=black,scale=\s] (1) at (.5,0){};
\node[circle, fill=black,scale=\s] (0) at (-.5,0){};
\draw[thick] (0) -- (1);
\node (lab) at (1,0){\footnotesize\textcolor{red}{-1}};
\node (labi) at (-1,0){\footnotesize\textcolor{white}{0}};
\end{tikzpicture}};
\node (2) at (0.5*\x,1*\y){\begin{tikzpicture}
\node[circle, fill=black,scale=\s] (1) at (.5,0){};
\node[circle, fill=black,scale=\s] (0) at (-.4,0){};
\node (lab) at (1,0){\footnotesize\textcolor{red}{-1}};
\node (labi) at (-1,0){\footnotesize\textcolor{white}{0}};
\end{tikzpicture}};
\node (1) at (0*\x,0.25*\y){\begin{tikzpicture}
\node[circle, fill=black,scale=\s] (0) at (0,0){};
\node (lab) at (.5,0){\footnotesize\textcolor{red}{1}};
\node (labi) at (-.5,0){\footnotesize\textcolor{white}{0}};
\end{tikzpicture}};
\draw[color=Gray] (H) -- (C4) -- (P3) -- (P2) --(1) -- (2) -- (21) -- (F) -- (H) -- (P4) -- (21) -- (P2);
\draw[color=Gray] (2) -- (P3) -- (F);
\draw[color=Gray] (P3) -- (P4);
\end{tikzpicture}
\caption{The interval $[K_1,H]$ in $\mathcal{G}$, where $H$ is the house graph, with $\mu(K_1,X)$ in red.}\label{fig:Hint}
\end{figure}

The \emph{dual} $P^*$ of a poset $P$ is the poset with the same elements of $P$ and the partial order is reversed,
 that is $a\le_{P^*} b$ if and only if $a\ge_P b$.  It is well known that $\mu_{P^*}(b,a)=\mu_P(a,b)$.
 We use $\mu$ to denote $\mu_{\mathcal{G}}$ and $\mu_c$ to denote $\mu_{\mathcal{G}^c}$,
 and $\mu^*$ and $\mu^*_c$ to denote M\"obius function on their duals.
 See \cref{fig:Hint} for an example of the M\"obius function on $\mathcal{G}$.

First we look
 at a condition on $G$ which results in $\mu(H,G)=0$ for all $|H|<|G|-1$. An \emph{automorphism}
 on a graph $G$ is a map $f$ from $V(G)$ onto itself such that $(u,v)$ is an edge if and only if~$(f(u),f(v))$ is an 
 edge. A graph is \emph{locally finite} if every vertex has a a finite number of neighbours.

\begin{defn}
A graph is \emph{vertex transitive} if for every pair of vertices $v_1$ and~$v_2$ there is an 
 automorphism that maps $v_1$ to $v_2$.
\end{defn}

The \emph{coatoms} of an interval $[H,G]$ are the maximal elements of $(H,G)$.
The following theorem allows us to consider what the coatoms of $[H,G]$ are
 when $G$ is vertex transitive.

\begin{thm}\label{thm:Tho87}\cite{Tho87}
Let $G$ be a locally finite graph without isolated vertices. Then $G$ is vertex transitive if and
 only if deleting any vertex gives the same graph up to isomorphism.
\end{thm}

If we restrict to only finite graphs then we can remove the condition on isolated vertices, as we now show.

\begin{cor}\label{cor:Tho87}
Let $G$ be a finite graph, then $G$ is vertex transitive if and only if deleting any vertex gives
 the same graph up to isomorphism.
\begin{proof}
If $G$ is finite and has no isolated vertices, then the result follows by \cref{thm:Tho87}. If~$G$ is
 finite and vertex transitive, then it must be regular, that is, all vertices have the same degree. So 
 if there is an isolated vertex all vertices must be isolated, so $G$ is the graph with no edges, which implies 
 all vertex deleted subgraphs are isomorphic. The converse follows by a similar argument.
\end{proof}
\end{cor}

\cref{cor:Tho87} implies that if $G$ is a vertex transitive graph then it contains, up to isomorphism, exactly one graph of order one
 less. Moreover, it is well known that if an interval $[x,y]$ has a single coatom, then $\mu(x,y)=0$, which immediately implies the following result.

\begin{prop}\label{prop:vtrans}
If a graph $G$ is vertex transitive, then $\mu(H,G)=0$ for any~${|H|<|G|-1}$.
\end{prop} 

A graph is \emph{simple} if it has no loops nor multiple edges. The \emph{complement} $\bar{G}$ of a simple graph $G$ has the same vertex set as $G$ and $(a,b)$ is an edge in
 $\bar{G}$ if and only if it is not an edge in $G$. Applying the complement operation to an interval of simple 
 graphs does not change the M\"obius function.

\begin{lem}
If $H$ and $G$ are simple graphs, then $[H,G]$ is isomorphic to $[\bar{H},\bar{G}]$, so $\mu(H,G)=\mu(\bar{H},\bar{G})$.
\begin{proof}
It is trivial to see that $g_1,\ldots,g_{|H|}$ is an occurrence of $H$ in $G$ if and only if it is an occurrence
 of $\bar{H}$ in $\bar{G}$. Therefore, $H\le G$ if and only if $\bar{H}\le\bar{G}$, which implies the result.
\end{proof}
\end{lem}

\subsection{Well known graphs}\label{sec:known}

In this subsection we consider the M\"obius function of intervals of some well known graphs. First note that the
 complete, empty and cycle graphs are vertex transitive, so by \cref{prop:vtrans} we get:
\begin{lem}\label{lem:KKC}
Consider any $n>0$ and graph $H$ with $|H|\not\in\{n-1,n\}$, then:
$$\mu(H,K_n)=\mu(H,\overline{K_n})=\mu(H,C_n)=0.$$
\end{lem}

Next we consider the intervals $[\emptyset,G]$, where $\emptyset$ is the null graph, that is, the graph with no vertices.
\begin{lem}\label{lem:emptyset}
If a graph $G$ has no loops, then:
$$\mu(\emptyset,G)=\begin{cases}
1,&\mbox{ if }|G|=0\\
-1,&\mbox{ if }|G|=1\\
0,&\mbox{ if }|G|>1
\end{cases}.$$
\begin{proof}
The cases $|G|\le 1$ follow trivially. Suppose $|G|>1$, then the interval~$(\emptyset,G)$ has a unique
 minimal element $K_1$, which implies $\mu(\emptyset,G)=0$.
\end{proof}
\end{lem}

If $G$ contains loops then \cref{lem:emptyset} does not hold as there are multiple minimal elements in~$(\emptyset,G)$,
 which are the graphs with a single vertex and $x$ loops, for some $x\ge 0$.

Next we consider bipartite graphs.

\begin{lem}\label{lem:K1K2}
If $G$ is a non-empty bipartite graph with $|G|>2$, then $$\mu(K_1,G)=-\mu(K_2,G).$$
\begin{proof}
Note that $K_2$ and $\overline{K_2}$ are contained in $G$ and in every element of $(K_2,G)$. The only  elements in
 $(K_1,G)\setminus(K_2,G)$ are $K_2$, $\overline{K_2}$ and possibly some larger empty graphs. Note that 
 $\mu(K_1,\overline{K_a})=0$, for any $a>2$, by \cref{lem:KKC}. Let~$Q$ be the poset obtained by removing all
 $\overline{K_a}$, with $a>2$, from $(K_1,G)$, then~$\mu_Q(K_1,G)=\mu(K_1,G)$. The poset $Q$ consists of $(K_2,G)\cup\{K_2,\overline{K_2}\}$,
 and every element of $(K_2,G)$ contains both~$K_2$ and $\overline{K_2}$, we consider the dual poset~$Q^*$ and let $\mathcal{I}=[G,\overline{K_2}]\cup K_2$, then
 $$\mu(K_1,G)=\mu_{Q^*}(G,K_1)=-\sum_{X\in \mathcal{I}}\mu_{Q^*}(G,X)=-\mu_{Q^*}(G,K_2)=-\mu(K_2,G).$$
\end{proof}
\end{lem}

\begin{rem}
We can also prove \cref{lem:K1K2} using a topological argument similar to that used in \cite[Theorem 4.2]{Smith14} to prove a result
 on the permutation pattern poset. To do so we show that $(K_1,G)$ is homotopically equivalent to $Q$ by the Quillen Fiber lemma, see \cite[Theorem~15.28]{Koz08}, and $Q$ is a
 suspension of $(K_2,G)$, which implies the result.
\end{rem}

Let $B_{a,b}$ be the complete bipartite graph with parts of size $a$ and $b$.
\begin{lem}\label{lem:bipart}
$$\mu(B_{a_1,a_2},B_{b_1,b_2})=\begin{cases}\displaystyle
1,&\mbox{ if }a_1=b_1\,\&\,a_2=b_2\\
-1,&\mbox{ if }(b_1-a_1)+(b_2-a_2)=1\\
1,&\mbox{ if }b_1-a_1=b_2-a_2=1\mbox{ and }a_1\not=a_2\\
0,&\mbox{ otherwise}
\end{cases}$$
\begin{proof}
The first two cases follow trivially as they are rank $0$ and $1$ intervals. If~$b_1-a_1=b_2-a_2=1$ the interval has rank $2$ and
 $(B_{a_1,a_2},B_{b_1,b_2})$ contains~$B_{a_1+1,a_2}$ and $B_{a_1,a_2+1}$, but if $a_1=a_2$ then these are isomorphic, 
 this implies the result for rank $2$ intervals. 

Suppose the interval has rank $r>2$. If $a_1=a_2$ or $b_1=b_2$, then there is a unique minimal element $B_{a_1+1,a_2}$ or maximal
 element $B_{b_1-1,b_2}$, respectively, which implies the result. Otherwise, by induction we can see that the only elements in 
 $[B_{a_1,a_2},B_{b_1,b_2})$ with non-zero M\"obius function are $B_{a_1,a_2}$, $B_{a_1+1,a_2}$, $B_{a_1,a_2+1}$ and $B_{a_1+1,a_2+1}$. Therefore,
 $\mu(B_{a_1,a_2},B_{b_1,b_2})
  =1-1-1+1=0.$
\end{proof}
\end{lem}

We can combine  \cref{lem:bipart} with \cref{lem:K1K2} to get  the following corollary.

\begin{cor}
For any $a,b$ with $a+b>2$:
$$\mu(K_1,B_{a,b})=0$$
\end{cor}

We can also consider the complete multipartite graphs $B_{b_{1}^{i_{1}},\ldots,b_{k}^{i_{k}}}$
 with $i_j$ partitions of size~$b_j$, for all $1\le j\le k$. 
\begin{lem}
Given any integers $a,n,k>0$ we have $\mu(H,B_{k^n})=0$, for all~${|H|<kn-1}$, and $\mu(B_{1^{a+n}},B_{1^a,k^n})=0$.
\begin{proof}
The graph $B_{k^n}$ is vertex transitive, which gives the first part of the statement. The only vertices that can be deleted
 from $B_{1^a,k^n}$ are those in the parts of size $k$, and deleting any vertex from any of those parts 
 gives the same graph up to isomorphism. So there is a unique maximal element.
\end{proof} 
\end{lem}

\section{Posets between graphs of disjoint paths}\label{sec:paths}

In this section we consider graphs which are a collection of disjoint paths, or equivalently acyclic graphs where every vertex has degree at most $2$.
 Recall that~$P_a$ is the path graph of order~$a$. Given a set $S$ of positive integers we define $\PG{S}$ to be the graph $\sqcup_{i\in S}P_i$, where $\sqcup$ denotes the disjoint union.
 We show that the M\"obius function from $\overline{K_n}$ to $n$ disjoint copies of $P_5$ is equal to the~$n$'th Catalan number $\mathcal{C}_n$. To show this we use discrete Morse theory, which we introduce
 in \cref{sec:sub2}. We also show that the M\"obius function from $\overline{K_2}$ to two disjoint copies of $P_n$ is equal to the $n$'th Fibonacci number $\mathcal{F}_n$.

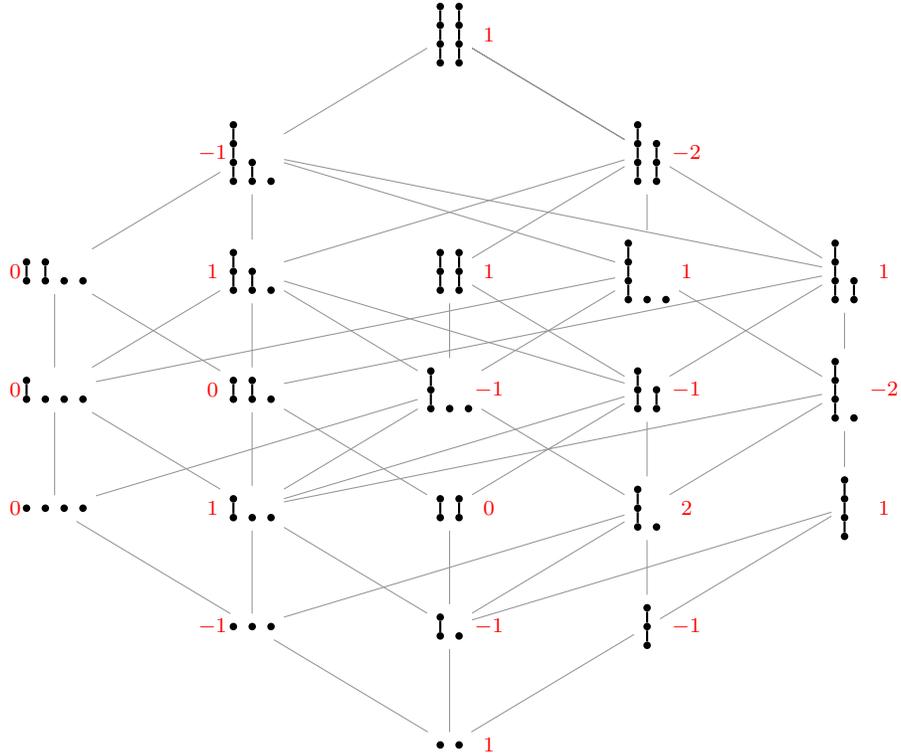
\begin{figure}\centering
\begin{tikzpicture}[scale=2.625]\def\y{0.6}\def\s{0.3}\def\sc{1}
\node (61) at (0,6*\y){\begin{tikzpicture}[scale=\sc]
\node[circle, fill=black,scale=\s] (1) at (0,0){};
\node[circle, fill=black,scale=\s] (2) at (0,0.25){};
\node[circle, fill=black,scale=\s] (3) at (0.25,0){};
\node[circle, fill=black,scale=\s] (4) at (0.25,0.25){};
\node[circle, fill=black,scale=\s] (5) at (0,0.5){};
\node[circle, fill=black,scale=\s] (6) at (0,0.75){};
\node[circle, fill=black,scale=\s] (7) at (0.25,0.5){};
\node[circle, fill=black,scale=\s] (8) at (0.25,0.75){};
\draw[thick] (1) -- (2) -- (5) -- (6);
\draw[thick] (3) -- (4) -- (7) -- (8);
\end{tikzpicture}};
\node (m61) at (0.2,6*\y){\footnotesize $\textcolor{red}{1}$};
\node (51) at (-1,5*\y){\begin{tikzpicture}[scale=\sc]
\node[circle, fill=black,scale=\s] (1) at (0,0){};
\node[circle, fill=black,scale=\s] (2) at (0,0.25){};
\node[circle, fill=black,scale=\s] (3) at (0.25,0){};
\node[circle, fill=black,scale=\s] (4) at (0.5,0){};
\node[circle, fill=black,scale=\s] (5) at (0.25,0.25){};
\node[circle, fill=black,scale=\s] (6) at (0,0.5){};
\node[circle, fill=black,scale=\s] (7) at (0,0.75){};
\draw[thick] (1) -- (2) -- (6) -- (7);
\draw[thick] (3) -- (5);
\end{tikzpicture}};
\node (m51) at (-1.2,5*\y){\footnotesize $\textcolor{red}{-1}$};
\node (52) at (1,5*\y){\begin{tikzpicture}[scale=\sc]
\node[circle, fill=black,scale=\s] (1) at (0,0){};
\node[circle, fill=black,scale=\s] (2) at (0,0.25){};
\node[circle, fill=black,scale=\s] (3) at (0.25,0){};
\node[circle, fill=black,scale=\s] (4) at (0.25,0.25){};
\node[circle, fill=black,scale=\s] (5) at (0,0.5){};
\node[circle, fill=black,scale=\s] (6) at (0,0.75){};
\node[circle, fill=black,scale=\s] (7) at (0.25,0.5){};
\draw[thick] (1) -- (2) -- (5) -- (6);
\draw[thick] (3) -- (4) -- (7);
\end{tikzpicture}};
\node (m52) at (1.2,5*\y){\footnotesize $\textcolor{red}{-2}$};
\node (41) at (-2,4*\y){\begin{tikzpicture}[scale=\sc]
\node[circle, fill=black,scale=\s] (1) at (0,0){};
\node[circle, fill=black,scale=\s] (3) at (0.25,0){};
\node[circle, fill=black,scale=\s] (2) at (0.5,0){};
\node[circle, fill=black,scale=\s] (4) at (0.75,0){};
\node[circle, fill=black,scale=\s] (5) at (0,0.25){};
\node[circle, fill=black,scale=\s] (6) at (0.25,0.25){};
\draw[thick] (1) -- (5);\draw[thick] (3) -- (6);
\end{tikzpicture}};
\node (m41) at (-2.2,4*\y){\footnotesize $\textcolor{red}{0}$};
\node (42) at (-1,4*\y){\begin{tikzpicture}[scale=\sc]
\node[circle, fill=black,scale=\s] (1) at (0,0){};
\node[circle, fill=black,scale=\s] (2) at (0,0.25){};
\node[circle, fill=black,scale=\s] (3) at (0.25,0){};
\node[circle, fill=black,scale=\s] (4) at (0.5,0){};
\node[circle, fill=black,scale=\s] (5) at (0.25,0.25){};
\node[circle, fill=black,scale=\s] (6) at (0,0.5){};
\draw[thick] (1) -- (2) -- (6);
\draw[thick] (3) -- (5);
\end{tikzpicture}};
\node (m42) at (-1.2,4*\y){\footnotesize $\textcolor{red}{1}$};
\node (43) at (0,4*\y){\begin{tikzpicture}[scale=\sc]
\node[circle, fill=black,scale=\s] (1) at (0,0){};
\node[circle, fill=black,scale=\s] (2) at (0,0.25){};
\node[circle, fill=black,scale=\s] (3) at (0.25,0){};
\node[circle, fill=black,scale=\s] (4) at (0.25,.25){};
\node[circle, fill=black,scale=\s] (5) at (0,0.5){};
\node[circle, fill=black,scale=\s] (6) at (0.25,0.5){};
\draw[thick] (1) -- (2) -- (5);
\draw[thick] (3) -- (4) -- (6);
\end{tikzpicture}};
\node (m43) at (0.2,4*\y){\footnotesize $\textcolor{red}{1}$};
\node (44) at (1,4*\y){\begin{tikzpicture}[scale=\sc]
\node[circle, fill=black,scale=\s] (1) at (0,0){};
\node[circle, fill=black,scale=\s] (2) at (0,0.25){};
\node[circle, fill=black,scale=\s] (3) at (0.25,0){};
\node[circle, fill=black,scale=\s] (4) at (0.5,0){};
\node[circle, fill=black,scale=\s] (5) at (0,0.5){};
\node[circle, fill=black,scale=\s] (6) at (0,0.75){};
\draw[thick] (1) -- (2) -- (5) -- (6);
\end{tikzpicture}};
\node (m44) at (1.2,4*\y){\footnotesize $\textcolor{red}{1}$};
\node (45) at (2,4*\y){\begin{tikzpicture}[scale=\sc]
\node[circle, fill=black,scale=\s] (1) at (0,0){};
\node[circle, fill=black,scale=\s] (2) at (0,0.25){};
\node[circle, fill=black,scale=\s] (3) at (0.25,0){};
\node[circle, fill=black,scale=\s] (4) at (0.25,0.25){};
\node[circle, fill=black,scale=\s] (5) at (0,0.5){};
\node[circle, fill=black,scale=\s] (6) at (0,0.75){};
\draw[thick] (1) -- (2) -- (5) -- (6);\draw[thick] (3) -- (4);
\end{tikzpicture}};
\node (m45) at (2.2,4*\y){\footnotesize $\textcolor{red}{1}$};
\node (31) at (-2,3*\y){\begin{tikzpicture}[scale=\sc]
\node[circle, fill=black,scale=\s] (1) at (0,0){};
\node[circle, fill=black,scale=\s] (3) at (0.25,0){};
\node[circle, fill=black,scale=\s] (2) at (0.5,0){};
\node[circle, fill=black,scale=\s] (4) at (0.75,0){};
\node[circle, fill=black,scale=\s] (5) at (0,0.25){};
\draw[thick] (1) -- (5);
\end{tikzpicture}};
\node (m31) at (-2.2,3*\y){\footnotesize $\textcolor{red}{0}$};
\node (32) at (-1,3*\y){\begin{tikzpicture}[scale=\sc]
\node[circle, fill=black,scale=\s] (1) at (0,0){};
\node[circle, fill=black,scale=\s] (2) at (0,0.25){};
\node[circle, fill=black,scale=\s] (3) at (0.25,0){};
\node[circle, fill=black,scale=\s] (4) at (0.5,0){};
\node[circle, fill=black,scale=\s] (5) at (0.25,0.25){};
\draw[thick] (1) -- (2);
\draw[thick] (3) -- (5);
\end{tikzpicture}};
\node (m32) at (-1.2,3*\y){\footnotesize $\textcolor{red}{0}$};
\node (33) at (0,3*\y){\begin{tikzpicture}[scale=\sc]
\node[circle, fill=black,scale=\s] (1) at (0,0){};
\node[circle, fill=black,scale=\s] (2) at (0,0.25){};
\node[circle, fill=black,scale=\s] (3) at (0.25,0){};
\node[circle, fill=black,scale=\s] (4) at (0.5,0){};
\node[circle, fill=black,scale=\s] (5) at (0,0.5){};
\draw[thick] (1) -- (2) -- (5);
\end{tikzpicture}};
\node (m33) at (0.2,3*\y){\footnotesize $\textcolor{red}{-1}$};
\node (34) at (1,3*\y){\begin{tikzpicture}[scale=\sc]
\node[circle, fill=black,scale=\s] (1) at (0,0){};
\node[circle, fill=black,scale=\s] (2) at (0,0.25){};
\node[circle, fill=black,scale=\s] (3) at (0.25,0){};
\node[circle, fill=black,scale=\s] (4) at (0.25,0.25){};
\node[circle, fill=black,scale=\s] (5) at (0,0.5){};
\draw[thick] (1) -- (2) -- (5);\draw[thick] (3) -- (4);
\end{tikzpicture}};
\node (m34) at (1.2,3*\y){\footnotesize $\textcolor{red}{-1}$};
\node (35) at (2,3*\y){\begin{tikzpicture}[scale=\sc]
\node[circle, fill=black,scale=\s] (1) at (0,0){};
\node[circle, fill=black,scale=\s] (2) at (0,0.25){};
\node[circle, fill=black,scale=\s] (3) at (0,0.5){};
\node[circle, fill=black,scale=\s] (4) at (0,0.75){};
\node[circle, fill=black,scale=\s] (5) at (0.25,0){};
\draw[thick] (1) -- (2) -- (3) -- (4);
\end{tikzpicture}};
\node (m35) at (2.2,3*\y){\footnotesize $\textcolor{red}{-2}$};
\node (21) at (-2,2*\y){\begin{tikzpicture}[scale=\sc]
\node[circle, fill=black,scale=\s] (1) at (0,0){};
\node[circle, fill=black,scale=\s] (3) at (0.25,0){};
\node[circle, fill=black,scale=\s] (2) at (0.5,0){};
\node[circle, fill=black,scale=\s] (4) at (0.75,0){};
\end{tikzpicture}};
\node (m21) at (-2.2,2*\y){\footnotesize $\textcolor{red}{0}$};
\node (22) at (-1,2*\y){\begin{tikzpicture}[scale=\sc]
\node[circle, fill=black,scale=\s] (1) at (0,0){};
\node[circle, fill=black,scale=\s] (2) at (0,0.25){};
\node[circle, fill=black,scale=\s] (3) at (0.25,0){};
\node[circle, fill=black,scale=\s] (4) at (0.5,0){};
\draw[thick] (1) -- (2);
\end{tikzpicture}};
\node (m22) at (-1.2,2*\y){\footnotesize $\textcolor{red}{1}$};
\node (23) at (0,2*\y){\begin{tikzpicture}[scale=\sc]
\node[circle, fill=black,scale=\s] (1) at (0,0){};
\node[circle, fill=black,scale=\s] (2) at (0,0.25){};
\node[circle, fill=black,scale=\s] (3) at (0.25,0){};
\node[circle, fill=black,scale=\s] (4) at (0.25,0.25){};
\draw[thick] (1) -- (2);\draw[thick] (3) -- (4);
\end{tikzpicture}};
\node (m23) at (0.2,2*\y){\footnotesize $\textcolor{red}{0}$};
\node (24) at (1,2*\y){\begin{tikzpicture}[scale=\sc]
\node[circle, fill=black,scale=\s] (1) at (0,0){};
\node[circle, fill=black,scale=\s] (2) at (0,.25){};
\node[circle, fill=black,scale=\s] (3) at (0,.5){};
\node[circle, fill=black,scale=\s] (4) at (0.25,0){};
\draw[thick] (1) -- (2) -- (3);
\end{tikzpicture}};
\node (m24) at (1.2,2*\y){\footnotesize $\textcolor{red}{2}$};
\node (25) at (2,2*\y){\begin{tikzpicture}[scale=\sc]
\node[circle, fill=black,scale=\s] (1) at (0,0){};
\node[circle, fill=black,scale=\s] (2) at (0,0.25){};
\node[circle, fill=black,scale=\s] (3) at (0,0.5){};
\node[circle, fill=black,scale=\s] (4) at (0,0.75){};
\draw[thick] (1) -- (2) -- (3) -- (4);
\end{tikzpicture}};
\node (m25) at (2.2,2*\y){\footnotesize $\textcolor{red}{1}$};
\node (11) at (-1,1*\y){\begin{tikzpicture}[scale=\sc]
\node[circle, fill=black,scale=\s] (1) at (0,0){};
\node[circle, fill=black,scale=\s] (3) at (.25,0){};
\node[circle, fill=black,scale=\s] (2) at (.5,0){};
\end{tikzpicture}};
\node (m11) at (-1.2,1*\y){\footnotesize $\textcolor{red}{-1}$};
\node (12) at (0,1*\y){\begin{tikzpicture}[scale=\sc]
\node[circle, fill=black,scale=\s] (1) at (0,0){};
\node[circle, fill=black,scale=\s] (2) at (0,0.25){};
\node[circle, fill=black,scale=\s] (3) at (0.25,0){};
\draw[thick] (1) -- (2);
\end{tikzpicture}};
\node (m12) at (0.2,1*\y){\footnotesize $\textcolor{red}{-1}$};
\node (13) at (1,1*\y){\begin{tikzpicture}[scale=\sc]
\node[circle, fill=black,scale=\s] (1) at (0,0){};
\node[circle, fill=black,scale=\s] (2) at (0,0.25){};
\node[circle, fill=black,scale=\s] (3) at (0,.5){};
\draw[thick] (1) -- (2) -- (3);
\end{tikzpicture}};
\node (m13) at (1.2,1*\y){\footnotesize $\textcolor{red}{-1}$};
\node (01) at (0,0*\y){\begin{tikzpicture}[scale=\sc]
\node[circle, fill=black,scale=\s] (1) at (0,0){};
\node[circle, fill=black,scale=\s] (3) at (.25,0){};
\end{tikzpicture}};
\node (m01) at (.2,0*\y){\footnotesize $\textcolor{red}{1}$};
\draw[color=Gray] (01) -- (11) -- (21) -- (31) -- (41) -- (51) -- (61) -- (52) -- (42) -- (32) -- (22) -- (12) -- (01) -- (13) -- (25) -- (35) -- (45) -- (52) -- (61);
\draw[color=Gray] (11) -- (22) -- (31) -- (42) -- (51) -- (44) -- (31);
\draw[color=Gray] (11) -- (24) -- (33) -- (42) -- (34) -- (43) -- (52) -- (44) -- (35) -- (24) -- (12) -- (25);
\draw[color=Gray] (12) -- (23) -- (32) -- (41);
\draw[color=Gray] (13) -- (24) -- (34) -- (45) -- (51);
\draw[color=Gray] (21) -- (33) -- (43);
\draw[color=Gray] (22) -- (33) -- (44);
\draw[color=Gray] (35) -- (22) -- (34) -- (23);
\draw[color=Gray] (32) -- (45);
\end{tikzpicture}
\caption{The interval $[11_p,44_p]$, with $\mu(11_p,X)$ in red.}
\label{fig:11-44}
\end{figure}

For ease of notation we often write $\PG{S}$ as $\PG{\alpha}$, where $\alpha$ is the word consisting of the letters of $S$ in decreasing order.
 To avoid confusion we use upper and lower case letters to denote the sets and words, respectively. 
 We also use $\PN{x}{n}$ to represent the graph made of $n$ disjoint copies of~$P_x$, and~$\PNN{x}{n}{y}{m}=\PN{x}{n}\sqcup\PN{y}{m}$.
 For example, the graph $\{4,4\}_P$, also written $\PG{{4,4}}$ or $\PN{4}{2}$, is two disjoint paths of length~$4$, see \cref{fig:11-44} for the interval $[\PG{{1,1}},\PG{{4,4}}]$.
 Let $\mathcal{P}$ be the subposet of $\mathcal{G}$ consisting of all graphs of disjoint paths, and let~$\mathcal{P}_{\le i}$ be the subposet
 of $\mathcal{P}$ of graphs where all paths have max length $i$.

In \cref{tab:2} we display the M\"obius function $\mu(\PN{1}{n},\PN{x}{n})$, for $x+n\le10$. We can see that for $n=2$ we get the Fibonacci numbers $\mathcal{F}_{x-2}$
 and for $x=5$ we get the Catalan numbers $\mathcal{C}_n$. In this section we prove the following result which deals with the first five columns and 
 first two rows of \cref{tab:2}.
 
  \begin{table}\centering
\begin{tabular}{c | *{10}{c}}
\diagbox{$n$}{$x$} & 1 &  2 & 3 &  4 & 5  &  6  &  7  &  8 & 9\\\hline							
1	 		   & 1 & -1 & 1 & -1 & 1  & -1  &  1  & -1 & 1\\
2	 		   & 1 &  0 & 1 &  1 & 2  &  3  &  5  &  8 &  \\
3	 		   & 1 &  0 & 1 & -1 & 5  & -14 &  47 &    &  \\	
4	 		   & 1 &  0 & 1 &  1 & 14 &  81 & 	  &    &  \\		
5	 		   & 1 &  0 & 1 & -1 & 42 &     & 	  &    &  \\ 			
6	 		   & 1 &  0 & 1 &  1 & 	  & 	&  	  &    &  \\
7	 		   & 1 &  0 & 1 &    &    &     & 	  &    &  \\					
8	 		   & 1 &  0 &   &    &    &     & 	  &    &  \\						
9	 		   & 1 &    &   &    &    &     &     &    &  \\
\end{tabular}
\caption{The values of $\mu(\PN{1}{n},\PN{x}{n})$, for $n+x\le10$.}\label{tab:2}
\end{table}

\begin{thm}\label{thm:fivecols}
Let $\mathcal{C}_i$ and $\mathcal{F}_i$ be the $i$'th Catalan and  Fibonacci numbers, respectively. 
 We can compute $\mu(\PN{1}{n},\PN{x}{n})$ for the following cases:
 $$\mu(\PN{1}{n},\PN{x}{n})=\begin{cases}
  (-1)^{x+n},&\mbox{ if } n=1\mbox{ or }x=4\\
  1,&\mbox{ if } x=1\mbox{ or }x=3\\
  0,&\mbox{ if } x=2\,\,\,\,\,\,\&\,\,\,\,\,\,n>1,\\
  F_{x-2},&\mbox{ if } n=2,\\
  \mathcal{C}_n,&\mbox{ if } x=5.\\
 \end{cases}$$
\end{thm}

An immediate corollary of \cref{thm:fivecols} is:
\begin{cor}
The M\"obius function is unbounded on $\mathcal{G}$.
\end{cor}

In the following two subsections we prove \cref{thm:fivecols} by
 \cref{cor:path,prop:fibcase,lem:5,lem:34,lem:12}. In \cref{sec:sub1} we apply an inductive argument and in 
 \cref{sec:sub2} we use discrete Morse theory.

\subsection{The cases \texorpdfstring{$n=1$}{n=1} and \texorpdfstring{$n=2$}{n=2}}\label{sec:sub1}
We begin with the case $n=1$, that is, the interval $[\PG{1},\PG{x}]$ from the singleton graph to the path graph of order $x$.
 To prove this case we prove a more general result on intervals $[\PG{m},\PG{x}]$ between any two path graphs, from which the result for $n=1$ immediately follows by setting $m=1$.
\begin{lem}\label{lem:path}
For any $0<m\le x$ we have $\mu(\PG{m},\PG{x})=(-1)^{x-m}$.
\begin{proof}
We consider the dual poset and prove $\ms(\PG{x},\PG{m})=(-1)^{x-m}$. Fixing~$x$, we proceed by induction decreasing the value of $m$. Clearly $\ms(\PG{x},\PG{x})=1$ and assume 
 $\ms(\PG{x},\PG{\ell})=(-1)^{x-m-1}$, for $\ell=m+1$. In the dual poset the graph $\PG{m}$ covers exactly two graphs $\PG{\ell}$ and $\PG{{m,1}}$, which
 is the graph obtained by adding an isolated vertex to $\PG{m}$. It is straightforward to see $Y\les\PG{{m,1}}$ for every $Y\in(\PG{x},\PG{m})^*\setminus \PG{\ell}$. So:
 \begin{align*}\ms(\PG{x},\PG{m})=-\left(\sum_{Y\in[\PG{x},\PG{{m,1}}]^*}\ms(\PG{x},Y)+\ms(\PG{x},\PG{\ell})\right)&=-\ms(\PG{x},\PG{\ell})\\&=-(-1)^{x-m-1}.\end{align*}
\end{proof}
\end{lem}
\begin{cor}\label{cor:path}
For any $x\ge1$ we have $\mu(\PG{1},\PG{x})=(-1)^{x-1}$.
\end{cor}

Next we prove the Fibonacci case $n=2$. For ease of notation we denote~$\mu(\PG{{1,1}},X)$ by $\tmu(X)$. First we need some Lemmas.
 In the proofs of the next four results we apply a similar argument each time, which we outline here. Given an element in an interval $z\in[x,y]$,
 then by the definition of the M\"obius function we know that $\sum_{w\in[x,z]}\mu(x,w)=0$, and when computing $\mu(x,y)$ we can ignore any $v\in[x,y]$ with $\mu(x,v)=0$.
 Therefore, $$\mu(x,y)=-\sum_{v\in L^z_{x,y}}\mu(x,v),$$ where we define
 $$L^z_{x,y}:=\{v\in[x,y)\,|\,v\not\le z\,\&\,\mu(x,v)\not=0\}.$$
 So our approach is to find an appropriate $z$, ignore everything in $[x,z]$ and look at all the remaining elements~$v$, ignoring
 any for which $\mu(x,v)=0$.
\begin{lem}\label{lem:mostzero}
If $A$ is a set with $|A|=3$ and $1\not\in A$ or $|A|>3$, then~${\tmu(\PG{A})=0}$.
\begin{proof}
The minimal cases for $A$ are $\{1,1,1,1\}$ or $\{2,2,2\}$.
 It is easy to see that~$\tmu(\PG{{1,1,1,1}})=\tmu(\PG{{2,2,2}})=0$.
 Assume the claim is true for any element smaller than $\PG{A}$ and consider $I=(\PG{{1,1}},\PG{{A}})$.
 Let $a$ and $b$ be the two largest elements of $A$.
 By induction for every element $X\in I$ with $\tmu(X)\not=0$ we have~$X\le\PG{{a,b,1}}$.
 Therefore, $L^{\PG{{a,b,1}}}_I$ is empty, which implies $\tmu(\PG{A})=0$.
\end{proof}
\end{lem}

\begin{lem}\label{lem:addone}
For any $a\ge b>1$:
$$\tmu(\PG{{a,b,1}})=-\tmu(\PG{{a,b}}).$$
\begin{proof}
\cref{fig:11-44} shows that $\tmu(\PG{{2,2,1}})=-\tmu(\PG{{2,2}})=0$.
Assume the claim is true for any $\PG{{i,j}}<\PG{{a,b}}$ and consider the interval $I=(\PG{{1,1}},\PG{{a,b,1}})$.
 Let ${Z=\PGs{{a,b-1,1}}}$, by \cref{lem:mostzero} we know that
$$L^Z_I=\{\PG{{i,b,1}}\,|\,b\le i< a\}\cup\{\PG{{j,b}}\,|\,b\le j\le a\}.$$ The induction hypothesis implies that 
 $\tmu(\PG{{j,b}})+\tmu(\PG{{j,b,1}})=0$, for $b\le j< a$. Therefore, ${\tmu(\PG{{a,b,1}})=-\tmu(\PG{{a,b}})}$.
\end{proof}
\end{lem}

\begin{lem}\label{lem:incone}
For any $a-1>b>1$:
$$\tmu(\PG{{a,b}})=-\tmu(\PGs{{a-1,b}}).$$
\begin{proof}
We can see in \cref{fig:11-44} that $\tmu(\PG{{4,2}})=\tmu(\PG{{3,2}})$.
 Assume the claim is true for $\PG{{i,j}}<\PG{{a,b}}$ and
 let $I=(\PG{{1,1}},\PG{{a,b}})$. Let $Z=\PGs{{a,b-2,1}}$,
 by \cref{lem:mostzero} the elements of $L^Z_I$ are of the form $\PG{{i,j}}$ or $\PG{{i,j,1}}$ and must have $i,j\ge b-1$. So:
 \begin{align*}L^Z_I=&\{\PGs{{a-1,b}},\PGs{{a-1,b-1}},\PGs{{a,b-1}}\}\\&\cup\{\PGs{{i,j}},\PGs{{i,j,1}}\,|\,b-1\le i\le a-2,\, j\in\{b,b-1\}\}\end{align*}
 \cref{lem:addone} and the induction hypothesis imply:
 \begin{align*}\tmu(\PG{{i,j}})+\tmu(\PG{{i,j,1}})&=0,\\\tmu(\PGs{{a-1,b-1}})+\tmu(\PGs{{a,b-1}})&=0.\end{align*}
 Therefore, $\tmu(\PG{{a,b}})=-\tmu(\PGs{{a-1,b}}).$
\end{proof}
\end{lem}

\begin{prop}\label{prop:fibcase}
For any $x>1$:
$$\tmu(\PGs{{x,x-1}})=-\mathcal{F}_{x-1}\hskip 20pt\text{ and }\hskip 20pt\tmu(\PG{{x,x}})=\mathcal{F}_{x-2}.$$
\begin{proof}
\cref{fig:11-44} shows that $\tmu(\PG{{2,2}})=\mathcal{F}_0$, $\tmu(\PG{{3,3}})=\mathcal{F}_1$,
 $\tmu(\PG{{2,1}})=-\mathcal{F}_1$ and $\tmu(\PG{{3,2}})=-\mathcal{F}_2$. 
 Assume the claim is true for any $i$ with~$3<i<x$. Let $I=(\PG{{1,1}},\PGs{{x,x-1}})$ and $Z=\PGs{{x,x-3,1}}$, also let $J=(\PG{{1,1}},\PG{{x,x}})$ and~${K=\PGs{{x,x-2,1}}}$.
 \cref{lem:mostzero} allows us to compute the sets:
 \begin{align*}L^{Z}_I&=\{\PGs{{x,x-2}}, \PGs{{x-1,x-1}}, \PGs{{x-1,x-2}}, \PGs{{x-2,x-2}},\\ &\hskip 150pt\PGs{{x-1,x-2,1}}, \PGs{{x-2,x-2,1}}\},\\[10pt]
 L^K_J&=\{\PGs{{x,x-1}},\PGs{{x-1,x-1}}\}.\end{align*}
 By \cref{lem:addone}, \cref{lem:incone} and the induction hypothesis we get 
 \begin{align*}&\tmu(\PGs{{x-1,x-2}})+\tmu(\PGs{{x-1,x-2,1}})=0,\\
 &\tmu(\PGs{{x-2,x-2}})+\tmu(\PGs{{x-2,x-2,1}})=0,\\
 &\tmu(\PGs{{x,x-2}})=-\tmu(\PGs{{x-1,x-2}})=\mathcal{F}_{x-2},\\
 &\tmu(\PGs{{x-1,x-1}})=\mathcal{F}_{x-3}.\end{align*} Therefore, 
 \begin{align}\label{eq:xx1}\tmu(\PGs{{x,x-1}})=-\tmu(\PGs{{x,x-2}})-\tmu(\PGs{{x-1,x-1}})&=-\mathcal{F}_{x-2}-\mathcal{F}_{x-3}\nonumber\\&=-\mathcal{F}_{x-1}.\end{align}
 Moreover, \cref{eq:xx1} then allows us to compute
 \begin{align*}\tmu(\PG{{x,x}})&=-\tmu(\PGs{{x,x-1}})-\tmu(\PGs{{x-1,x-1}})=\mathcal{F}_{x-1}-\mathcal{F}_{x-3}=\mathcal{F}_{x-2}.\end{align*}
\end{proof}
\end{prop}
Combining \cref{lem:incone} and \cref{prop:fibcase} gives the following corollary:
\begin{cor}
For any $a>b>1$:
$$\tmu(\PG{{a,b}})=(-1)^{a+b}\mathcal{F}_{b}.$$
\end{cor}

\subsection{The cases \texorpdfstring{$1\le x\le 5$}{1 ≤ x ≤ 5}}\label{sec:sub2}

To prove the results
 in this subsection we use discrete Morse theory, first developed in \cite{For95}. In particular we use the methods introduced in \cite{BabHer05} for applying discrete Morse
 theory to the order complex of a poset. We give a brief introduction to the necessary theory, but for further details we refer the reader to \cite[Section 4]{SagVat06}.
 In this subsection we consider the dual poset exclusively, so we let $\le$ denote the dual partial order and drop the $*$ superscript. 
 Also, we use $A$, $B$ and $C$ to denote chains and $a_i$, $b_i$ and $c_i$ to denote the $i$'th element of these chains, respectively.
 
 \mbox{}

Consider an interval $[x,y]$. We say two chains diverge at index $i$ if they agree up to $i$ but the $i+1$'th elements are different.
 An ordering $\lhd$ of the maximal chains of $(x,y)$ is a \emph{PL-ordering} if given two chains $A\lhd B$ which diverge at $i$ and two chains $A'$ and $B'$
 which agree to index $i+1$ with $A$ and $B$, respectively, then $A'\lhd B'$.
 
Let $C\setminus \hat{C}$ be the chain obtained by deleting
 $\hat{C}$ from $C$. Consider a PL-ordering $\lhd$ of the maximal chains of $(x,y)$, along with a maximal chain $C$ and subchain $(c_i,c_{i+k+1})_C:=(c_{i+1}<c_{i+2}<\cdots<c_{i+k})$.
 If $C\setminus (c_i,c_{i+k+1})_C$ is a subchain of a chain $B\lhd C$,
 then we say $(c_i,c_{i+k+1})$ is a \emph{skipped interval} of $C$, of size $k$. A skipped interval is a \emph{minimal skipped interval (MSI)} if it does not
 strictly contain another skipped interval.
 
Let $I_1,I_2,\ldots,I_t$ be the MSIs of a chain $C$ in increasing order of the index of their first element.
 Set $J_1=I_1$, and then $I'_k=I_k\setminus J_1$, for all $k>1$, remove any~$I'_k$ which is no longer minimal,
 and set the first remaining one as $J_2$. Repeat this process until there are no non-empty modified MSIs remaining, and denote the set of $J_1,J_2,\ldots,J_{k'}$ by $\mathcal{J}(C)$.
 A chain $C$ is \emph{critical} if every element of $C$ appears in $\mathcal{J}(C)$.
 
Given a PL-ordering and the sets $\mathcal{J}(C)$ we can compute the M\"obius function~$\mu(x,y)$ which equals the
 number of critical chains in any PL-ordering, with the sign given by the rank, see \cite[Proposition 3.1]{BabHer05}. Another useful result is \cite[Proposition 4.2]{BabHer05}
 which states that a PL-ordering is a shelling order if
 for every chain $C$ the set $\mathcal{J}(C)$ only contains skipped intervals of size $1$. 
 
 \mbox{}

To apply discrete Morse theory we need to introduce a PL-ordering for the chains of $\pf$. Consider two graphs $\PG{{\alpha}}\lessdot\PG{{\beta}}$ in $\pa$,
 to get from $\PG{{\alpha}}$ to $\PG{{\beta}}$
 we must delete a vertex of $\PG{{\alpha}}$ from a path of length $u\ge1$.
 Deleting this point will take the component $\PG{{u}}$ to $\PG{{w^1,w^2}}$, where $w^1+w^2=u-1$ and $w^1\ge w^2\ge0$.
 Define the \emph{operation} between $\PG{{\alpha}}$ and $\PG{{\beta}}$ as the pair $u$ and $w:=\{w^1,w^2\}$ which we
 call the \emph{domain} and \emph{image} of the operation, respectively. We denote the operation by $\oper{u}{w^1w^2}$, or
 as $\oper{u}{w_1}$ when $w_2=0$.
 
We can represent any maximal chain by the operations applied at each step. Given a maximal chain $C=c_1<c_2<\dots<c_n$, of an interval of $\pa$, define the \emph{operation chain}
 $\Lambda(C)=(\lc{C}{1},\lc{C}{2},\ldots,\lc{C}{n-1})$, where $\lc{C}{i}$ is the operation that takes $c_i$ to $c_{i+1}$. 
 For example,  $C=\PG{{5,1}}<\PG{{2,2,1}}<\PG{{2,1,1}}<\PG{{2,1}}$ has the operation chain $\Lambda(C)=(\fitwtw,\twon,\onze)$.
 
 We can define a PL-ordering of the maximal chains of any interval of $\pa$ by defining an ordering $\prec$ on the set of all possible 
 operations. To do this take two maximal chains $B$
 and $C$ which diverge at $i$ and order $B$ before $C$, denoted~$B\lhd C$, if and only if $\lc{B}{i}\prec\lc{C}{i}$.
 There are $9$ possible operations in $\pf$, which we give the following order:
 \begin{equation}\fitwtw\prec\fotwon\prec\thtw\prec\fifo\prec\fithon\prec\foth\prec\thonon\prec\onze\prec\twon.\end{equation}

In \cref{tab:1} we represent all 81 possible pairs of operations $\lc{C}{i},\lc{C}{i+1}$, where each cell is split into three parts. The top part of the cell contains
 the chain $c_i<c_{i+1}<c_{i+2}$, where $c_i$ contains exactly the domains of $\lc{C}{i}$ and $\lc{C}{i+1}$, and~${c_{i+1}<c_{i+2}}$ is obtained by applying $\lc{C}{i}$ and $\lc{C}{i+1}$ to $c_i$.
 If the domain of $\lc{C}{i+1}$ is contained in the image 
 of $\lc{C}{i}$, then there is a second chain where $c_i$ is exactly the domain of $\lc{C}{i}$. The middle part of the cell
 contains all elements~${b_1,\ldots,b_k}$ that can be obtained from $c_i$ by operations ordered before $\lc{C}{i}$. The bottom part of the cell contains 
 the operations $O_1,\ldots,O_k$, where $O_j$ maps~$b_j$ to $c'_{i+2}$, if no such operation exists then $O_j=\emptyset$.

\begin{table}[ht]
\centering
\resizebox{\textwidth}{!}{
\begin{tabular}{|c|c|c|c|c|c|c|c|c|c|}
\arrayrulecolor{black}\cline{1-10}
\diagbox{$\lc{C}{i}$}{\raisebox{-5pt}{$\lc{C}{i+1}$}} 	&{\raisebox{-2.5pt}{\huge$\fitwtw$}}&\raisebox{-2.5pt}{\huge$\fotwon$}   &\raisebox{-2.5pt}{\huge$\thtw$} 		 &\raisebox{-2.5pt}{\huge$\fifo$} 		&\raisebox{-2.5pt}{\huge$\fithon$}		&\raisebox{-2.5pt}{\huge$\foth$} 		&\raisebox{-2.5pt}{\huge$\thonon$} 	  &\raisebox{-2.5pt}{\huge$\onze$} 		&\raisebox{-2.5pt}{\huge$\twon$} 		 \\ \arrayrulecolor{black}\cline{1-10}
\lev         											&$55<522<2222$   				    &$54<422<2221$         				 &$53<322<222$           				 &$55<522<422$          				&$55<522<3221$           				&$54<422<322$           			    &$53<322<2211$         				  &$51<221<22$         					&\tl{$52<222<221$}{$5<22<21$}            \\ \arrayrulecolor{gray}\cdashline{2-10}
{\huge$\fitwtw$}                             			&     				   			    &         				 			 &         				     			 &         								&           				  			&          				   				&        				  			  &        								&\tl{}{}			 					 \\ \arrayrulecolor{gray}\cdashline{2-10}
\lev                                         			&          			   		   		&         				 			 &           				 			 &          							&           				  			&           				    		&         				  			  &       								&\tl{}{}		     					 \\ \arrayrulecolor{black}\cline{1-10}
\lev						           		 			&$54<521<2221$           		    &$44<421<2211$         				 &$43<321<221$           				 &$54<521<421$          				&$54<521<3211$           				&$44<421<321$           				&$43<321<2111$         				  &\tl{$41<211<21$}{$4<21<2$}           &\tl{$42<221<211$}{$4<21<11$}            \\ \arrayrulecolor{gray}\cdashline{2-10}
{\huge$\fotwon$}                             			&$422$           		   		    &						     		 &				             			 &$422$	                        		&$422$						          	&			           				    &			         				  &\tl{}{}        						&\tl{}{}          						 \\ \arrayrulecolor{gray}\cdashline{2-10}
\lev                                         			&$\fotwon$    					    &					         		 &          					 		 &$\twon$						        &$\emptyset$				            &			          					&							          &\tl{}{}        						&\tl{}{}           						 \\ \arrayrulecolor{black}\cline{1-10}
\lev             							 			&$53<52<222$           		   		&$43<42<221$        				 &$33<32<22$	           				 &$53<52<42$          					&$53<52<321$		           			&$43<42<32$           					&$33<32<211$         				  &$31<21<2$         					&\tl{$32<22<21$}{$3<2<1$}                \\ \arrayrulecolor{gray}\cdashline{2-10}
{\huge$\thtw$}                               			&$322$           				   	&$321$ 				         		 &				          				 &$322$       						    &$322$							  	  	&$321$           						&			         				  &			         					&\tl{}{}           						 \\ \arrayrulecolor{gray}\cdashline{2-10}
\lev                                         			&$\thtw$           			   		&$\thtw$					         &				           				 &$\emptyset$				            &$\twon$					          	&$\onze$				            	&			        				  &			         					&\tl{}{}           						 \\ \arrayrulecolor{black}\cline{1-10}
\lev             							 			&$55<54<422$           		   		&\tl{$54<44<421$}{$5<4<21$}        	 &$53<43<42$	           				 &$55<54<44$          				    &$55<54<431$           				  	&\tl{$54<44<43$}{$5<4<3$}               &$53<43<411$    			          &$51<41<4$         					&$52<42<41$         					 \\ \arrayrulecolor{gray}\cdashline{2-10}
{\huge$\fifo$}                               			&$522$				               	&\tl{$422,521$}{$22$}				 &$322,52$					             &$522$					            	&$522$		           			      	&\tl{$422,521$}{$22$}           		&$322,52$         	  				  &$221$				            	&$222$					             	 \\ \arrayrulecolor{gray}\cdashline{2-10}
\lev                                         			&$\fifo$          	   			   	&\tl{$\twon,\fifo$}{$\twon$}		 &$\emptyset,\fifo$           			 &$\emptyset$			        		&$\emptyset$					    	&\tl{$\emptyset,\emptyset$}{$\emptyset$}&$\emptyset,\emptyset$				  &$\emptyset$         					&$\emptyset$					 		 \\ \arrayrulecolor{black}\cline{1-10}
\lev           								 			&$55<531<3221$         		   		&$54<431<3211$				         &\tl{$53<331<321$}{$5<31<21$}			 &$55<531<431$          				&$55<531<3311$				          	&$54<431<331$           				&\tl{$53<331<3111$}{$5<31<111$}       &\tl{$51<311<31$}{$5<31<3$}           &$52<321<311$         					 \\ \arrayrulecolor{gray}\cdashline{2-10}
{\huge$\fithon$}                             			&$522,54$           			   	&$422,521,44$		                 &\tl{$322,52,43$}{$22,4$}       	 	 &$522,54$						        &$522,54$				         		&$422,521,44$           				&\tl{$322,52,43$}{$22,4$}         	  &\tl{$221,41$}{$22,4$}         		&$222,42$         						 \\ \arrayrulecolor{gray}\cdashline{2-10}
\lev                                         			&$\fithon,\emptyset$              	&$\emptyset,\fithon,\emptyset$		 &\tl{$\twon,\fithon,\fotwon$}{$\twon,\fotwon$}&$\emptyset,\fithon$   			&$\emptyset,\emptyset$					&$\emptyset,\emptyset,\emptyset$		&\tl{$\emptyset,\emptyset,\emptyset$}{$\emptyset,\emptyset$}&\tl{$\emptyset,\foth$}{$\emptyset,\foth$}&$\emptyset,\emptyset$   \\ \arrayrulecolor{black}\cline{1-10}
\lev             							 			&$54<53<322$           		   		&$44<43<321$         	 			 &\tl{$43<33<32$}{$4<3<2$}            	 &$54<53<43$           					&$54<53<331$				            &$44<43<33$	            				&\tl{$43<33<311$}{$4<3<11$}			  &$41<31<3$          					&$42<32<31$          					 \\ \arrayrulecolor{gray}\cdashline{2-10}
{\huge$\foth$}                               			&$422,521,44,431$		           	&$421$			 		 			 &\tl{$321,42$}{$21$}           		 &$422,44,521,431$				        &$422,521,44,431$			            &$421$           						&\tl{$321,42$}{$21$}      			  &$211$         						&$221$         							 \\ \arrayrulecolor{gray}\cdashline{2-10}
\lev                                         			&$\foth,\emptyset,\emptyset,\emptyset$&$\foth$			     			 &\tl{$\onze,\foth$}{$\onze$}  			 &$\emptyset,\foth,\emptyset,\onze$     &$\emptyset,\emptyset,\emptyset,\foth$	&$\emptyset$           					&\tl{$\twon,\emptyset$}{$\twon$} 	  &$\emptyset$					        &$\emptyset$					         \\ \arrayrulecolor{black}\cline{1-10}
\lev           								 			&$53<511<2211$           		   	&$43<411<2111$         				 &$33<311<211$           				 &$53<511<411$          			    &$53<511<3111$           				&$43<411<311$           			    &$33<311<1111$         				  &\tl{$31<111<11$}{$3<11<1$}           &$32<211<111$         					 \\ \arrayrulecolor{gray}\cdashline{2-10}
{\huge$\thonon$}                             			&$322,52,43,331$           	   		&$321,42,33$         				 &$32$				         			 &$322,52,43,331$				        &$322,52,43,331$						&$321,42,33$           					&$32$				          		  &\tl{$21$}{$2$}         				&$22$        						     \\ \arrayrulecolor{gray}\cdashline{2-10}
\lev                                         			&$\thonon,\emptyset,\emptyset,\emptyset$&$\thonon,\emptyset,\emptyset$   &$\thonon$								 &$\emptyset,\emptyset,\thonon,\emptyset$&$\emptyset,\emptyset,\emptyset,\thonon$&$\twon,\emptyset,\thonon$          	&$\emptyset$            			  &\tl{$\twon$}{$\twon$}         		&$\emptyset$         				     \\ \arrayrulecolor{black}\cline{1-10}
\lev             						  	 			&$51<5<22$           			   	&$41<4<21$        				     &$31<3<2$           					 & $51<5<4$         			        &$51<5<31$           					&$41<4<3$           					&$31<3<11$         					  &$11<1<0$         					&$21<2<1$         						 \\ \arrayrulecolor{gray}\cdashline{2-10}
{\huge$\onze$}                               			&$221,41,311$           		   	&$211,31$         					 &$21,111$					         	 &$221,41,311$ 					        &$221,41,311$           				&$211,31$					            &$21,111$         					  &			         					&								         \\ \arrayrulecolor{gray}\cdashline{2-10}
\lev                                         			&$\onze,\emptyset,\emptyset$      	&$\onze,\thtw$				         &$\onze,\emptyset$	             		 &$\emptyset,\onze,\emptyset$           &$\emptyset,\foth,\onze$				&$\emptyset,\onze$    	                &$\twon,\onze$			          	  &								        &							             \\ \arrayrulecolor{black}\cline{1-10}
\lev             							 			&$52<51<221$           		   		&$42<41<211$         				 &$32<31<21$           			     	 &$52<51<41$          					&$52<51<311$           				  	&$42<41<31$           					&$32<31<111$         				  &\tl{$21<11<1$}{$2<1<0$}              &$22<21<11$         					 \\ \arrayrulecolor{gray}\cdashline{2-10}
{\huge$\twon$}                               			&$222,42,321$                     	&$221,32$				             &$22,211$ 					         	 &$222,42,321$ 					        &$222,42,321$						    &$221,32$						        &$22,211$					          &\tl{$2$}{\mbox{}}        	    	&										 \\ \arrayrulecolor{gray}\cdashline{2-10}
\lev                                         			&$\twon,\fotwon,\thtw$            	&$\twon,\thonon$				     &$\twon,\onze$				             &$\emptyset,\twon,\emptyset$           &$\emptyset,\emptyset,\twon$            &$\emptyset,\twon$			            &$\emptyset,\twon$			          &\tl{$\twon$}{\mbox{}}            	&				             			 \\ \arrayrulecolor{black}\cline{1-10}
\end{tabular}}
\caption{
The table of all pairs $(\lc{C}{i},\lc{C}{i+1})$
where each cell is split into three parts. The top part is $c_i<c_{i+1}<c_{i+2}$. The middle part
 contains all elements $b_1,\ldots,b_k$ that can be obtained from $c_i$ by operations ordered before $\lc{C}{i}$.
 The bottom part contains 
 the operations~$O_1,\ldots,O_k$, where $O_j$ maps $b_j$ to $c'_{i+2}$, if no such operation exists then $O_j=\emptyset$.
 So $(c_i,c_{i+2})$ is an MSI if and only if the cell $(\lc{C}{i},\lc{C}{i+1})$ has at least one operation $O\not=\emptyset$
 in the bottom part of the cell. We omit the subscript $p$ in the table for ease of notation.
}\label{tab:1}
\end{table}

\cref{tab:1} allows us to see all MSIs of size $1$, but there are larger MSIs, such as $(\PG{{5,3}},\PG{{3,2}})$ which is an MSI of size $2$ of the chain $\PG{{5,3}}<\PG{{5,2}}<\PG{{4,2}}<\PG{{3,2}}$. 
 Therefore, this PL-ordering is not a shelling order. In fact, there is no ordering of the operations which induces a shelling order in this way, that is, for every ordering we can find a MSI of size greater than one.
 We have computationally verified this, the smallest interval needed to show it is $[\PG{{5,4,3}},\PG{{1,1,1}}]$\footnote{This code is available upon request.}.
 However, this does not imply that the intervals of $\pf$ are not shellable. In fact it is likely these intervals are all shellable, but a different approach is needed to prove this.
 
To show the M\"obius function results we need to count the critical chains, so first we establish some notation and present some results on what can be an MSI.
 We say that $C=c_i<c_{i+1}<\dots<c_{i+k+1}$ contains the operation $\oper{u}{w}$ if it appears in $\Lambda(C)$, and $C$ begins or ends with $\oper{u}{w}$ if $\lc{C}{i}=\oper{u}{w}$ or $\lc{C}{i+k}=\oper{u}{w}$,
 respectively. Similarly, we say a skipped interval $(c_i,c_j)$ contains, begins or ends with an operation if the subchain $c_{i+1}<\dots<c_{j-1}$ has that property.
 We denote by $C\oper{u}{w}$ the number of times $\oper{u}{w}$ appears in $\Lambda(C)$. 
 We can construct a chain $A$ by making changes to $\Lambda(C)$, but we must ensure that the new chain is a valid operation chain.
 An operation chain $\Lambda(A)$ is invalid if $\lc{A}{i}=\oper{u}{w}$ and there is no $u$ in $a_i$. Note that if $I$ is a skipped
 interval of $C$ then $I$ is a skipped interval of any chain $D$ which contains $C$, so when looking to see if
 $I=(c_i,c_j)$ appears as an MSI of a chain in $\pf$ it suffices to consider chains from $c_i$ to $c_j$.

\begin{lem}\label{no:order}
If $I=(c_i,c_{i+k+1})$ is a skipped interval of the chain $C=c_i<c_{i+1}<\dots<c_{i+k+1}$ in $\pf$, with $k>1$, and $B\lhd C$ is a chain from $c_i$ to $c_{i+k+1}$,
 then $I$ is not minimal if any of the following hold:
\begin{enumerate}
\item $C$ contains a pair $\lc{C}{y}\prec\lc{C}{x}$, where $x<y$ and $(\lc{C}{y},\lc{C}{x})\not=(\onze,\twon)$,   \label{no:1}
\item Either of $C$ or $B$ contain both $\thonon$ and $\onze$,\label{no:2}
\item $C$ contains both $\fithon$ and $\onze$\label{no:2a}
\item $B$ contains $\twon$ and $1\in c_{i+k+1}$,\label{no:3}
\item $B$ contains $\lc{C}{i}$.\label{no:4}
\end{enumerate}
\begin{proof}
If there is a pair $x<y$ with $\lc{C}{y}\prec\lc{C}{x}$ and $(\lc{C}{y},\lc{C}{x})\not=(\onze,\twon)$, then there is a pair $\lc{C}{\ell+1}\prec\lc{C}{\ell}$, with $x\le \ell<y$ and $(\lc{C}{\ell},\lc{C}{\ell+1})\not=(\onze,\twon)$.
 However, \cref{tab:1} shows that for any such $\lc{C}{\ell+1}\prec\lc{C}{\ell}$, that is, everything south east of the diagonal except $(\onze,\twon)$, we get an MSI $(c_\ell,c_{\ell+2})$. 
 So $I$ is not minimal, giving a contradiction.
 
Suppose $C$ contains $\thonon$ and $\onze$, with $\lc{C}{x}=\thonon$
 and $\lc{C}{y}=\onze$ the leftmost occurrence of each.
 If $y<x$, then $I$ is not minimal by \eqref{no:1}. If $x<y$, then we can create a new chain from $c_x$ to $c_{y+1}$ with the operation chain
 $\Lambda(A)=\thtw,\lc{C}{x+1},\lc{C}{x+2},\dots,\lc{C}{y-1},\twon$.
 We can see that $\Lambda(A)$ is a valid operation chain because the domain of $\lc{C}{x}$ and $\lc{A}{x}$ are the same, and
 their images contain the domains of $\lc{C}{y}$ and $\lc{A}{y}$, respectively.
 So, $a_{y+1}=c_{y+1}$ and $A$ is lexicographically less than $c_x<\dots<c_y$, thus $(c_x,c_y)$ is a skipped interval contained in $I$, so $I$ is not minimal.
 
Suppose $B$ contains $\thonon$ and $\onze$. If $\onze$ is before $\thonon$ in $\Lambda(B)$, then create a new operation chain $\Lambda(A)$ by moving $\onze$ to immediately after $\thonon$.
 This is a valid operation chain as no operation acts on the image of $\onze$. We can then apply the same argument to $\Lambda(A)$ as that used above, for the case $C$ contains both $\thonon$ and $\onze$,
 to see that $I$ is not minimal.
 
Suppose $C$ contains $\fithon$ and $\onze$. We can apply an analogous argument to~\eqref{no:2} where we create a new operation chain $\Lambda(A)$ from $\Lambda(C)$ by replacing 
 $\fithon$ and $\onze$ with $\fifo$ and $\foth$, respectively.
 
Suppose $B$ contains $\twon$. Let $\Lambda(A)$ be obtained from $\Lambda(B)$ by moving the~$\lc{B}{t}=\twon$ to the end. If this creates an invalid chain 
 then it must be caused by $\lc{B}{k}=\onze$, with $k>t$. So also move the $\onze$ to the penultimate position of~$\Lambda(A)$, which is valid 
 as $c_{i+k+1}$ contains a $1$ so there must be an operation~$\lc{B}{v}$ whose image contains a $1$, with $v>k$.
 We can see that $A\lhd B\lhd C$, because we move $\twon$ to the right of an operation ordered before it.
 Therefore,~$A$ implies that $(c_i,c_{i+k})$ is a skipped interval, contradicting the claim that $I$ is minimal.

Suppose $\Lambda(B)$ contains $\lc{C}{i}$. Let $\Lambda(A)$ be obtained from $\Lambda(B)$ by moving $\lc{C}{i}$ to the start. Clearly $b_{i}=c_i$ contains the domain of $\lc{C}{i}$
 so this is a valid operation chain. Therefore, 
 $(c_{i+1},c_{i+k+1})$ is a skipped interval, so $I$ is not minimal. 
\end{proof}
\end{lem}

\begin{lem}\label{lem:notwon}
There is no MSI in any chain of $\pf$ that ends in $\twon$.
\begin{proof}
Let $I$, $B$ and $C$ be as defined in \cref{no:order}.
 As $c_{i+k+1}=b_{i+k+1}$ we know that the number of times each letter $j$ appears must be the same in both, and the number of components must be the same in both.
 We can check this is the case by looking in $B$ and $C$ at the operations that have $j$ in their image or domain, and the operations that create or remove components.

For a contradiction assume that $\lc{C}{i+k}=\twon$. So we know that $1\in c_{i+k+1}$, thus \cref{no:order}\eqref{no:3} implies that $B$ does not contain $\twon$. There must be an operation in $C$
 whose image contains $2$, otherwise the number of $2$'s in $c_{i+k+1}$ and $b_{i+k+1}$ will differ.
 So $C$ must contain at least one of $\fitwtw,\fotwon$ or $\thtw$. In \cref{clm:1,clm:2,clm:3} we show that $C$ cannot contain
 any of these, so we get a contradiction.

\begin{clm}\label{clm:1}
The chain $C$ cannot contain $\fitwtw$.
\begin{proof}
Suppose $C$ contains $\fitwtw$, then the leftmost occurrence is either at the start or is preceded by a different operation.
 If $\lc{C}{i}=\fitwtw$, then $\lc{B}{i}=\fitwtw$ otherwise $C\lhd B$, so $I$ is not minimal by \cref{no:order}\eqref{no:4}.
 If $\fitwtw$ is preceded by a different operation $O$, then we must have $\fitwtw\prec O$, so $I$ is not minimal by \cref{no:order}\eqref{no:1}.
\end{proof}
\end{clm}

\begin{clm}\label{clm:2}
The chain $C$ cannot contain $\fotwon$.
\begin{proof}
Suppose $C$ contains $\fotwon$, and by \cref{clm:1} we can assume $C$ does not contain $\fitwtw$. So $\lc{C}{1}=\fotwon$ because otherwise there is a pair $(\alpha,\fotwon)$,
 with $\fotwon\prec\alpha$, which implies $I$ is not minimal by \cref{no:order}\eqref{no:1}.
 Moreover, $B$ cannot contain $\fotwon$ by \cref{no:order}\eqref{no:3}, so $\lc{B}{1}=\fitwtw$ as $B\lhd C$. Therefore, $C$ must contain some operation whose domain is $5$, so either $\fifo$ or $\fithon$, otherwise
 the number of $5$'s in $c_{i+k+1}$ and $b_{i+k+1}$ will differ.
 
If $C$ contains $\fifo$, then the leftmost $\fifo$ must be 
 preceded by $\thtw$, otherwise~$I$ is not minimal by \cref{tab:1}. So $\Lambda(C)$ must begin $\fotwon^u,\thtw^v,\fifo$,
 where~$u,v>0$ denote multiple consecutive operations. Let $\Lambda(A)$ be the chain obtained from~$\Lambda(C)$ by
 replacing the leftmost $\fotwon$ with $\fitwtw$ and $\fifo$ with 
 $\twon$. It is easy to verify that $\Lambda(A)$ a valid operation chain, and we can see that in both case we remove $5$ and add $2$ and $1$.
 So $A\lhd C$ and $a_{u+v+1}=c_{u+v+1}$, which implies $(c_i,c_{u+v+1})$ is a skipped interval contained in $I$, so $I$ is not minimal
 
If $C$ contains $\fithon$, and does not contain $\fifo$, then $\Lambda(C)$ must begin with $\fotwon^u,\fithon$. 
 By \cref{no:order}\eqref{no:2a} we know that $C$ cannot contain $\onze$ and by \cref{no:order}\eqref{no:2} 
 we know that $B$ cannot contain both $\thonon$ and $\onze$. So we consider the two cases where $B$ does or does not contain $\thonon$.

If $B$ does not contain $\thonon$, then counting the number of $1$'s added and the
 number of additional components created in $B$ and $C$, gives the following equalities:
\begin{align}
C\twon+\,2C\thonon+\,C\fotwon+\,C\fithon&=B\fithon-\,B\onze\label{num:1}\\
C\thonon+\,C\fotwon+\,C\fithon&=B\fithon-\,B\onze+\,B\fitwtw\label{num:2}\end{align}
Subtracting \eqref{num:2} from \eqref{num:1} implies $C\thonon+\,C\twon=-\,B\fitwtw$, which we know is not
 possible as all values are non-negative and $B\fitwtw>0$.

If $B$ contains $\thonon$, so does not contain $\onze$, then counting the components added, the number of $1$'s added,
 the number of $5$'s removed and the number of~$2$'s added in $B$ and $C$ implies:
\begin{align}C\thonon+\,C\fotwon+\,C\fithon&=B\fithon+\,B\thonon+\,B\fitwtw\label{num:3}\\
C\twon+\,2C\thonon+\,C\fotwon+\,C\fithon&=B\fithon+\,2B\thonon\label{num:4}\\
C\fithon&=B\fithon+\,B\fitwtw+\,B\fifo\label{num:5}\\
C\fotwon-\,C\twon&=B\thtw+\,2B\fitwtw\label{num:6}
\end{align}
Combining the above equations, by \eqref{num:6}+\eqref{num:4}+\eqref{num:5}-2\eqref{num:3}, gives $$B\fitwtw + \,B\thtw + \,B\fifo = 0,$$
 which again is impossible as $B\fitwtw>0$ and all values are non-negative.

Therefore, we cannot have $\fithon$ in $C$. So there is no operation in $C$ whose domain is $5$,
 which implies $c_{i+k+1}\not=b_{i+k+1}$, hence $C$ cannot start with $\fotwon$.
\end{proof}
\end{clm}

\begin{clm}\label{clm:3}
The chain $C$ cannot contain $\thtw$.
\begin{proof}
Assume $C$ contains $\thtw$ and does not contain $\fitwtw$ or $\fotwon$. So $C$ must begin with~$\thtw$ by \cref{no:order}\eqref{no:1}.
 First we show that $C$ cannot contain $\fithon$. If $\fithon$is in $\Lambda(C)$, then $\Lambda(C)$ must begin $\thtw^u,\fifo^v,\fithon$.
  Let $\Lambda(A)$ be obtained from $\Lambda(C)$ by replacing the leftmost occurrences of~$\thtw,\fifo$ and $\fithon$
  with~$\fitwtw,\twon$ and $\fifo$, respectively. It is straightforward to see that this is a valid operation chain. Moreover, $a_{i+u+v+1}=c_{i+u+v+1}$ and $A$ is a
  lexicographically less chain, thus $(c_i,c_{i+u+v+1})$ is a skipped interval contained in $I$, so $I$ is not minimal.
  
We can also see that
  there cannot be a $\foth$ in $C$. If there is then let $\Lambda(A)$ be obtained from $\Lambda(C)$ by changing $\thtw$ to $\fotwon$ and $\foth$ to $\onze$.
  So $A$ is lexicographically less, thus $I$ is not minimal.

By \cref{no:order} we can assume $B$ does not contain $\thtw$ nor both $\thonon$ and~$\onze$.
 If there is no $\thonon$ in $B$ then we must have at least as many $3$'s in $b_{i+k+1}$ as in $b_i$. But this is not possible, as $a_i$ and $a_{i+k+1}$ have a different number of $3$'s
 because we remove a $3$ from $c_i$ and have no
 way to create more as~$C$ does not contain $\fithon$ nor $\foth$. So there must be a $\thonon$ in $B$, which means there is no $\onze$. This implies there must be a $\thonon$ in $C$ as 
 $b_{i+k+1}$ has more components than $b_i$ so we need some operation in $C$ that creates components, and $\thonon$ is the only remaining option. 
 So $C$ can contain $\thtw,\fifo,\thonon$ and $\twon$. So counting the number of components added and $3$'s removed from $B$ and $C$, we get:
 \begin{align}C\thonon&=B\thonon+\,B\fithon+\,B\fitwtw+\,B\fotwon,\label{num:7}\\C\thtw+\,C\thonon&=B\thonon-\,B\fithon-\,B\foth.\label{num:8}\end{align}
  Subtracting \eqref{num:7} from \eqref{num:8}, and rearranging, we get:
 \begin{align}B\fitwtw+\,B\fotwon+\,2B\fithon+\,B\foth+\,C\thtw=0,\end{align}
 which is impossible, as we know that $C\thtw>0$ and all values are non-negative.
 Therefore, $C$ cannot begin $\thtw$.
\end{proof}
\end{clm}
\end{proof}
\end{lem}
%

\begin{lem}\label{lem:nofifo}
If an MSI of a chain in $\pf$ begins with $\fifo$, then it is of size~$1$.
\begin{proof}
Let $I$, $B$ and $C$ be as defined in \cref{no:order} and suppose $C$ begins with $\fifo$. 
 The only way to create a $4$ is by the operation $\fifo$ and by \cref{no:order}\eqref{no:4}, we 
 know that $B$ does not contain $\lc{C}{i}=\fifo$. Therefore, the number of $4$'s in~$b_{i+k+1}$ is weakly less than in $b_i$.
 So $C$ must contain some operation whose domain is $4$. By \cref{no:order}\eqref{no:1} this cannot be $\fotwon$, so must be $\foth$.
 Therefore,~$C$ cannot contain $\thonon$, because $(\foth,\thonon)$ is an MSI, and $C$ does not contain $\thtw$ because $\thtw\prec\fifo$, so we cannot remove $3$'s.
 
The possible operations with domain $5$ in $C$ are $\fithon$ and $\fifo$. Moreover, for every $\fifo$ we then must have a $\foth$.
 This implies that for every $5$ removed from~$c_i$ there is an extra $3$ in $c_{i+k+1}$. However, $B$ begins with $\fitwtw$ so this cannot be the case in $B$.
 So no valid $B$ exists, giving a contradiction to $I$ being a skipped interval
\end{proof}
\end{lem}

We now have the necessary information on the MSIs of $\pf$ to compute the critical chains, and prove 
 the cases $1\le x\le 5$ of \cref{thm:fivecols}.

\begin{prop}\label{lem:5}
For any $n\ge1$:
$$\mu(\PN{1}{n},\PN{5}{n})=\mathcal{C}_n.$$
\begin{proof}
To show the result we need to compute all critical chains.
 By \cref{lem:notwon} we know that a critical chain $C$ cannot contain $\twon$, hence $C$ also cannot contain~$\fitwtw,\fotwon$ nor $\thtw$, as $\PN{1}{n}$ does not contain any $2$'s.
  Also note that $C$ cannot contain $\fifo$, because either $\fifo$ is in the middle of an MSI or at the end of one MSI and beginning of another.
  The second case is not possible because \cref{lem:nofifo} and \cref{tab:1} imply
  in any MSI beginning with $\fifo$ it is followed by one of $\fitwtw$, $\fotwon$ or $\thtw$. The first case is not possible because \cref{tab:1} implies
  if~$\fifo$ is preceded by anything other than $\fitwtw$, $\fotwon$ or $\thtw$ then it creates an MSI which ends with $\fifo$. Therefore, $C$ also cannot contain $\foth$ as there is never a $4$ in any element of the chain. 
 So the only operations in $C$ are $\onze$,~$\fithon$, and~$\thonon$. 
 
We can see that $C\fithon=n$, $C\thonon=n$ and $C\onze=2n$, in order to get from~$\PN{5}{n}$ to $\PN{1}{n}$.
 Note that $\onze$ cannot be preceded by another~$\onze$, as this cannot create an MSI. So $\onze$ is preceded by $\fithon$ or $\thonon$
 and as $C\fithon+\,C\thonon=C\onze$ we know that every $\fithon$ and $\thonon$ is followed by $\onze$, so we have no choice in the placement of $\onze$.
 Secondly note that we can only apply $\thonon$ if there is a $3$ for it to be applied to, which implies that to the left of every $\thonon$ there is at least one more $\fithon$ than $\thonon$.
 \cref{tab:1} implies that any chain of this form is critical.
 Therefore, if we let $\fithon$ denote an north east step and $\thonon$ denote a south east step, then the critical
 chains are in bijection with the Dyck paths of length $2n$. Finally, the rank is equal to 
 $4n$, so the sign is given by~$(-1)^{4n}=1$.
\end{proof}
\end{prop}

\begin{prop}\label{lem:34}
For any non-negative integers $x+y=n$:
 $$\mu(\PN{1}{n},\PNN{4}{x}{3}{y})=(-1)^{x}\dbinom{n}{y}.$$
\begin{proof}
By an analogous argument to that used in the proof of \cref{lem:5} we can see that a critical chain $C$ can only contain $\foth,\thonon$ and $\onze$. Moreover,~$C$ must have $C\foth=x$, $C\thonon=n$ and~$C\onze=n$. Note that $\onze$  must be preceded by $\thonon$, 
 and as $C\thonon=C\onze$ this implies every $\thonon$ must be followed by a $\onze$. So we have $n$ copies of the consecutive pair $\thonon,\onze$ and we can choose to insert
 $\foth$ before any $\thonon$ so there are $\binom{n}{y}$ such chains. \cref{tab:1} implies that any chain of this form is critical.
 Finally, the rank equals $3x+2y$ so the sign is given by $(-1)^{3a+2b}=(-1)^a$.
\end{proof}
\end{prop}

\begin{lem}\label{lem:12}
For any $n>0$ we have
$\mu(\PN{1}{n},\PN{1}{n})=1$ and $$\mu(\PN{1}{n},\PN{2}{n})=\begin{cases}
-1,&\mbox{ if }n=1\\
0,&\mbox{ if }n>1
\end{cases}.$$
\begin{proof}
It is trivial to see that $\mu(\PN{1}{n},\PN{1}{n})=1$. The interval $(\PN{1}{n},\PN{2}{n})$ is a chain of the graphs $\PNN{2}{x}{1}{y}$, with $x+y=n$, which implies the result.
\end{proof}
\end{lem}

\cref{lem:5,lem:34,lem:12} imply cases $1\le x\le 5$ of \cref{thm:fivecols}, so along with \cref{sec:sub1} we have finished the proof of \cref{thm:fivecols}.

We finish this section with two conjectures. Computational evidence indicates that the value of $\mu(\PN{1}{n},\PNN{5}{x}{4}{y})$ is given by the Schr\"oder numbers, see sequence A088617 in \cite{OEIS}.
\begin{conj}
For any non-negative integers $x+y=n$:
 $$\mu(\PN{1}{n},\PNN{5}{x}{4}{y})=(-1)^{y}T(n,x),$$
 where $T(n,x)$ are the Schr\"oder numbers.
\end{conj}

There is a similar poset to $[\PN{1}{n},\PN{5}{n}]$ for which the M\"obius function also appears to be the Catalan numbers.
 Let $H^n$ be $n$ disjoint copies of the house graph, which is the graph at the top of \cref{fig:Hint}.

\begin{conj}
For any $n>0$:
$$\mu(\PN{1}{n},H^n)=\mu(\PN{1}{n},\PN{5}{n})=\mathcal{C}_n.$$
\end{conj}

\section{Further Questions}\label{sec:conj}
In this section we finish with some open problems and conjectures relating to the poset of graphs.
 A natural question to ask is  for what proportion of intervals is the M\"obius function non-zero. 
 If we take all intervals $[H,G]$ of simple graphs with $|G|\le n$ with $n=4,5,6,7$ the proportion 
 of intervals with $\mu(H,G)=0$ are approximately: $16.7\%$, $19.1\%$, $19.0\%$ and $14.2\%$, respectively.
 However this is a small sample and we do not know what happens for larger graphs.
\begin{que}
What proportion of intervals have non-zero M\"obius function?
\end{que}

\begin{figure}[ht]\centering$\begin{tikzpicture}
\node[circle, fill=black,scale=0.5] (1) at (0,3.5){};
\node[circle, fill=black,scale=0.5] (2) at (-.5,3){};
\node[circle, fill=black,scale=0.5] (3) at (.5,3){};
\node[circle, fill=black,scale=0.5] (4) at (-.5,2){};
\node[circle, fill=black,scale=0.5] (5) at (.5,2){};
\node[circle, fill=black,scale=0.5] (6) at (0,1.5){};
\node[circle, fill=black,scale=0.5] (7) at (0,1){};
\draw[thick] (1) -- (2) -- (4) -- (5) -- (3) -- (1);
\draw[thick] (2) -- (5);
\draw[thick] (4) -- (3);
\draw[thick] (4) -- (6) -- (5);
\draw[thick] (6) -- (7);
\end{tikzpicture}$
\caption{A graph $G$ for which $\mu(\overline{K_2},G)=1$.}\label{fig:nonalt}
\end{figure}
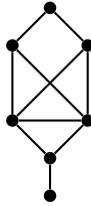

We say the M\"obius function is alternating if the sign is given by the rank, that is $\mu(H,G)$ has sign $(-1)^{|G|-|H|}$.
 This is not always true in $\mathcal{G}$. For example, let $G$ be the graph in \cref{fig:nonalt}, then $\mu(\bar{K}_2,G)=1$ 
 so the sign is not given by~$(-1)^{7-2}$.

If the sign of the M\"obius function is not determined by the parity of the rank, then the interval is not shellable. We know that an interval is not
 shellable if it contains a non-trivial disconnected subinterval. Are there any other obstructions to 
 shellability in $\mathcal{G}$ and $\mathcal{G}^c$
 and what proportion of intervals are shellable?

For many pattern posets it is conjectured that rank functions of all intervals are unimodal, but this is often difficult to prove.
\begin{conj}
All intervals of $\mathcal{G}$ and $\mathcal{G}^c$ are unimodal.
\end{conj}

Studying intervals of $\mathcal{G}$ has links to the graph reconstruction problem, which states that every graph $G$ is
 uniquely determined by the set of graphs obtained by deleting a single vertex from $G$. If the reconstruction conjecture 
 is true, it implies the coatoms of $[K_1,G]$ are unique for every $G$ with $|G|>2$, where the coatoms are the maximal 
 elements of the interior. But is this true if we change the bottom graph of the interval?
\begin{que}
Are the coatoms of $[H,G]$ unique for every pair $H<G$?
\end{que}

 \newcommand{\noop}[1]{}

\end{document}